# A new approach to inverse spectral theory, II. General real potentials and the connection to the spectral measure*

By Fritz Gesztesy and Barry Simon

## Abstract

We continue the study of the $A$-amplitude associated to a half-line Schrödinger operator, $-\frac{d^2}{dx^2} + q$ in $L^2((0,b))$, $b \leq \infty$. $A$ is related to the Weyl-Titchmarsh $m$-function via $m(-\kappa^2) = -\kappa - \int_0^a A(\alpha) e^{-2\alpha\kappa}\, d\alpha + O(e^{-(2a-\varepsilon)\kappa})$ for all $\varepsilon > 0$. We discuss five issues here. First, we extend the theory to general $q$ in $L^1((0,a))$ for all $a$, including $q$'s which are limit circle at infinity. Second, we prove the following relation between the $A$-amplitude and the spectral measure $\rho$: $A(\alpha) = -2 \int_{-\infty}^{\infty} \lambda^{-\frac{1}{2}} \sin(2\alpha\sqrt{\lambda})\, d\rho(\lambda)$ (since the integral is divergent, this formula has to be properly interpreted). Third, we provide a Laplace transform representation for $m$ without error term in the case $b < \infty$. Fourth, we discuss $m$-functions associated to other boundary conditions than the Dirichlet boundary conditions associated to the principal Weyl-Titchmarsh $m$-function. Finally, we discuss some examples where one can compute $A$ exactly.

## 1. Introduction

In this paper we will consider Schrödinger operators

$$(1.1) \qquad -\frac{d^2}{dx^2} + q$$

in $L^2((0,b))$ for $0 < b < \infty$ or $b = \infty$ and real-valued locally integrable $q$. There are essentially four distinct cases.

*This material is based upon work supported by the National Science Foundation under Grant No. DMS-9707661. The government has certain rights in this material.

1991 *Mathematics Subject Classification*. Primary: 34A55, 34B20; Secondary: 34L05, 47A10.

*Key words and phrases*. Inverse spectral theory, Weyl-Titchmarsh $m$-function, spectral measure.

*Case* 1.   $b < \infty$. We suppose $q \in L^1((0,b))$. We then pick $h \in \mathbb{R} \cup \{\infty\}$ and add the boundary condition at $b$

$$(1.2) \qquad u'(b_-) + hu(b_-) = 0,$$

where $h = \infty$ is shorthand for the Dirichlet boundary condition $u(b_-) = 0$.

For Cases 2–4, $b = \infty$ and

$$(1.3) \qquad \int_0^a |q(x)|\, dx < \infty \qquad \text{for all } a < \infty.$$

*Case* 2.   $q$ is "essentially" bounded from below in the sense that

$$(1.4) \qquad \sup_{a>0}\left(\int_a^{a+1} \max(-q(x),0)\, dx\right) < \infty.$$

Examples include $q(x) = c(x+1)^\beta$ for $c > 0$ and all $\beta \in \mathbb{R}$ or $q(x) = -c(x+1)^\beta$ for all $c > 0$ and $\beta \leq 0$.

*Case* 3.   (1.4) fails but (1.1) is limit point at $\infty$ (see [6, Ch. 9]; [33, Sect. X.1] for a discussion of limit point/limit circle), that is, for each $z \in \mathbb{C}_+ = \{z \in \mathbb{C} \mid \operatorname{Im}(z) > 0\}$,

$$(1.5) \qquad -u'' + qu = zu$$

has a unique solution, up to a multiplicative constant, which is $L^2$ at $\infty$. An example is $q(x) = -c(x+1)^\beta$ for $c > 0$ and $0 < \beta \leq 2$.

*Case* 4.   (1.1) is limit circle at infinity; that is, every solution of (1.5) is $L^2((0,\infty))$ at infinity if $z \in \mathbb{C}_+$. We then pick a boundary condition by picking a nonzero solution $u_0$ of (1.5) for $z = i$. Other functions $u$ satisfying the associated boundary condition at infinity then are supposed to satisfy

$$(1.6) \qquad \lim_{x \to \infty}[u_0(x)u'(x) - u_0'(x)u(x)] = 0.$$

Examples include $q(x) = -c(x+1)^\beta$ for $c > 0$ and $\beta > 2$.

The Weyl-Titchmarsh $m$-function, $m(z)$, is defined for $z \in \mathbb{C}_+$ as follows. Fix $z \in \mathbb{C}_+$. Let $u(x,z)$ be a nonzero solution of (1.5) which satisfies the boundary condition at $b$. In Case 1, that means $u$ satisfies (1.2); in Case 4, it satisfies (1.6); and in Cases 2–3, it satisfies $\int_R^\infty |u(x,z)|^2\, dx < \infty$ for some (and hence for all) $R \geq 0$. Then,

$$(1.7) \qquad m(z) = \frac{u'(0_+, z)}{u(0_+, z)}$$

and, more generally,

$$(1.8) \qquad m(z, x) = \frac{u'(x, z)}{u(x, z)}.$$

$m(z, x)$ satisfies the Riccati equation (with $m' = \frac{\partial m}{\partial x}$),

$$(1.9) \qquad m'(z, x) = q(x) - z - m(z, x)^2.$$

$m$ is an analytic function of $z$ for $z \in \mathbb{C}_+$, and moreover:

*Case 1.* $m$ is meromorphic in $\mathbb{C}$ with a discrete set $\lambda_1 < \lambda_2 < \cdots$ of poles on $\mathbb{R}$ (and none on $(-\infty, \lambda_1)$).

*Case 2.* For some $\beta \in \mathbb{R}$, $m$ has an analytic continuation to $\mathbb{C} \setminus [\beta, \infty)$ with $m$ real on $(-\infty, \beta)$.

*Case 3.* In general, $m$ cannot be continued beyond $\mathbb{C}_+$ (there exist $q$'s where $m$ has a dense set of polar singularities on $\mathbb{R}$).

*Case 4.* $m$ is meromorphic in $\mathbb{C}$ with a discrete set of poles (and zeros) on $\mathbb{R}$ with limit points at both $+\infty$ and $-\infty$.

Moreover,
$$\text{if } z \in \mathbb{C}_+ \text{ then } m(z, x) \in \mathbb{C}_+;$$

so $m$ satisfies a Herglotz representation theorem,

$$(1.10) \qquad m(z) = c + \int_{\mathbb{R}} \left[ \frac{1}{\lambda - z} - \frac{\lambda}{1 + \lambda^2} \right] d\rho(\lambda),$$

where $\rho$ is a positive measure called the spectral measure, which satisfies

$$(1.11) \qquad \int_{\mathbb{R}} \frac{d\rho(\lambda)}{1 + |\lambda|^2} < \infty,$$

$$(1.12) \qquad d\rho(\lambda) = \text{w-lim}_{\varepsilon \downarrow 0} \frac{1}{\pi} \text{Im}(m(\lambda + i\varepsilon)) \, d\lambda,$$

where w-lim is meant in the distributional sense.

All these properties of $m$ are well known (see, e.g. [23, Ch. 2]).

In (1.10), $c$ (which is equal to $\text{Re}(m(i))$) is determined by the result of Everitt [10] that for each $\varepsilon > 0$,

$$(1.13) \quad m(-\kappa^2) = -\kappa + o(1) \quad \text{as } |\kappa| \to \infty \text{ with } -\frac{\pi}{2} + \varepsilon < \arg(\kappa) < -\varepsilon < 0.$$

Atkinson [3] improved (1.13) to read,

$$(1.14) \qquad m(-\kappa^2) = -\kappa + \int_0^{a_0} q(\alpha) e^{-2\alpha\kappa} \, d\alpha + o(\kappa^{-1})$$

again as $|\kappa| \to \infty$ with $-\frac{\pi}{2} + \varepsilon < \arg(\kappa) < -\varepsilon < 0$ (actually, he allows $\arg(\kappa) \to 0$ as $|\kappa| \to \infty$ as long as $\text{Re}(\kappa) > 0$ and $\text{Im}(\kappa) > -\exp(-D|\kappa|)$ for suitable $D$). In (1.14), $a_0$ is any fixed $a_0 > 0$.

One of our main results in the present paper is to go way beyond the two leading orders in (1.14).

THEOREM 1.1. *There exists a function $A(\alpha)$ for $\alpha \in [0, b)$ so that $A \in L^1((0, a))$ for all $a < b$ and*

$$(1.15) \qquad m(-\kappa^2) = -\kappa - \int_0^a A(\alpha) e^{-2\alpha\kappa} \, d\alpha + \tilde{O}(e^{-2a\kappa})$$

*as $|\kappa| \to \infty$ with $-\frac{\pi}{2} + \varepsilon < \arg(\kappa) < -\varepsilon < 0$. Here we say $f = \tilde{O}(g)$ if $g \to 0$ and for all $\varepsilon > 0$, $(f/g)|g|^\varepsilon \to 0$ as $|\kappa| \to \infty$. Moreover, $A - q$ is continuous and*

$$(1.16) \qquad |(A - q)(\alpha)| \leq \left[ \int_0^\alpha |q(x)| \, dx \right]^2 \exp\left( \alpha \int_0^\alpha |q(x)| \, dx \right).$$

This result was proven in Cases 1 and 2 in [35]. Thus, one of our purposes here is to prove this result if one only assumes (1.3) (i.e., in Cases 3 and 4).

Actually, in [35], (1.15) was proven in Cases 1 and 2 for $\kappa$ real with $|\kappa| \to \infty$. Our proof under only (1.3) includes Case 2 in the general $\kappa$-region $\arg(\kappa) \in (-\frac{\pi}{2} + \varepsilon, -\varepsilon)$ and, as we will remark, the proof also holds in this region for Case 1.

*Remark.* At first sight, it may appear that Theorem 1.1 as we stated it does not imply the $\kappa$ real result of [35], but if the spectral measure $\rho$ of (1.10) has $\mathrm{supp}(\rho) \in [a, \infty)$ for some $a \in \mathbb{R}$, (1.15) extends to all $\kappa$ in $|\arg(\kappa)| < \frac{\pi}{2} - \varepsilon$, $|\kappa| \geq a + 1$. To see this, one notes by (1.10) that $m'(z)$ is bounded away from $[a, \infty)$ so one has the *a priori* bound $|m(z)| \leq C|z|$ in the region $\mathrm{Re}(z) < a - 1$. This bound and a Phragmén-Lindelöf argument let one extend (1.15) to the real $\kappa$ axis.

Here is a result from [35] which we will need:

THEOREM 1.2 ([35, Theorem 2.1]). *Let $q \in L^1((0, \infty))$. Then there exists a function $A(\alpha)$ on $(0, \infty)$ so that $A - q$ is continuous and satisfies (1.16) such that for $\mathrm{Re}(\kappa) > \frac{1}{2}\|q\|_1$,*

$$(1.17) \qquad m(-\kappa^2) = -\kappa - \int_0^\infty A(\alpha) e^{-2\alpha\kappa} \, d\alpha.$$

*Remark.* In [35], this is only stated for $\kappa$ real with $\kappa > \frac{1}{2}\|q\|_1$, but (1.16) implies that $|A(\alpha) - q(\alpha)| \leq \|q\|_1^2 \exp(\alpha\|q\|_1)$ so the right-hand side of (1.17) converges to an analytic function in $\mathrm{Re}(\kappa) > \frac{1}{2}\|q\|_1$. Since $m(z)$ is analytic in $\mathbb{C} \setminus [\alpha, \infty)$ for suitable $\alpha$, we have equality in $\{\kappa \in \mathbb{C} \mid \mathrm{Re}(\kappa) > \frac{1}{2}\|q\|_1\}$ by analyticity.

Theorem 1.1 in all cases follows from Theorem 1.2 and the following result which we will prove in Section 3.

THEOREM 1.3. *Let $q_1, q_2$ be potentials defined on $(0, b_j)$ with $b_j > a$ for $j = 1, 2$. Suppose that $q_1 = q_2$ on $[0, a]$. Then in the region $\arg(\kappa) \in (-\frac{\pi}{2} + \varepsilon, -\varepsilon)$, $|\kappa| \geq K_0$, we have that*

$$|m_1(-\kappa^2) - m_2(-\kappa^2)| \leq C_{\varepsilon,\delta} \exp(-2a\operatorname{Re}(\kappa)), \tag{1.18}$$

*where $C_{\varepsilon,\delta}$ depends only on $\varepsilon$, $\delta$, and $\sup_{0 \leq x \leq a}(\int_x^{x+\delta} |q_j(y)| \, dy)$, where $\delta > 0$ is any number so that $a + \delta \leq b_j$, $j = 1, 2$.*

*Remarks.* 1. An important consequence of Theorem 1.3 is that if $q_1(x) = q_2(x)$ for $x \in [0, a]$, then $A_1(\alpha) = A_2(\alpha)$ for $\alpha \in [0, a]$. Thus, $A(\alpha)$ is only a function of $q$ on $[0, \alpha]$. At the end of the introduction, we will note that $q(x)$ is only a function of $A$ on $[0, x]$.

2. This implies Theorem 1.1 by taking $q_1 = q$ and $q_2 = q\chi_{[0,a]}$ and using Theorem 1.2 on $q_2$.

3. Our proof implies (1.18) on a larger region than $\arg(\kappa) \in (-\frac{\pi}{2} + \varepsilon, -\varepsilon)$. Basically, we will need $\operatorname{Im}(\kappa) \geq -C_1 \exp(-C_2|\kappa|)$ if $\operatorname{Re}(\kappa) \to \infty$.

We will obtain Theorem 1.3 from the following pair of results.

THEOREM 1.4. *Let $q$ be defined on $(0, a+\delta)$ and $q \in L^1((0, a+\delta))$. Then in any region $\arg(\kappa) \in (-\frac{\pi}{2} + \varepsilon, -\varepsilon)$, $|\kappa| \geq K_0$, we have for all $x \in [0, a]$ that*

$$|m(-\kappa^2, x) + \kappa| \leq C_{\varepsilon,\delta}, \tag{1.19}$$

*where $C_{\varepsilon,\delta}$ depends only on $\varepsilon, \delta$ and $\sup_{0 \leq x \leq a}(\int_x^{x+\delta} |q(y)| \, dy)$.*

THEOREM 1.5. *Let $q_1 = q_2$ on $[0, a]$ and suppose $m_1$ and $m_2$ obey (1.19) for $x \in [0, a]$. Then in the same $\kappa$-region,*

$$|m_1(-\kappa^2) - m_2(-\kappa^2)| \leq 2C_{\varepsilon,\delta} \exp(-(\operatorname{Re}(\kappa))(2a - 2C_{\varepsilon,\delta})). \tag{1.20}$$

We will prove Theorem 1.5 in Section 2 using the Riccati equation and Theorem 1.4 in Section 3 by following ideas of Atkinson [3].

In Sections 5–9, we turn to the connection between the spectral measure $d\rho$ and the $A$-amplitude. Our basic formula says that

$$A(\alpha) = -2 \int_{-\infty}^{\infty} \lambda^{-\frac{1}{2}} \sin(2\alpha\sqrt{\lambda}) \, d\rho(\lambda). \tag{1.21}$$

In this formula, if $\rho$ gives nonzero weight to $(-\infty, 0]$, we interpret

$$\lambda^{-\frac{1}{2}} \sin(2\alpha\sqrt{\lambda}) = \begin{cases} 2\alpha & \text{if } \lambda = 0, \\ (-\lambda)^{-\frac{1}{2}} \sinh(2\alpha\sqrt{-\lambda}) & \text{if } \lambda < 0, \end{cases} \tag{1.22}$$

consistent with the fact that $\lambda^{-\frac{1}{2}} \sin(2\alpha\sqrt{\lambda})$ defined on $(0, \infty)$ extends to an entire function of $\lambda$.

The integral in (1.21) is not convergent. Indeed, the asymptotics (1.13) imply that $\int_0^R d\rho(\lambda) \sim \frac{2}{3\pi} R^{\frac{3}{2}}$ so (1.21) is never absolutely convergent. As we will see in Section 9, it is never even conditionally convergent in case $b < \infty$ (and also in many cases with $b = \infty$). So (1.21) has to be suitably interpreted.

In Sections 5–7, we prove (1.21) as a distributional relation, smeared in $\alpha$ on both sides by a function $f \in C_0^\infty((0, \infty))$. This holds for all $q$'s in Cases 1–4. In Section 8, we prove an Abelianized version of (1.21), *viz.*,

$$(1.23) \qquad A(\alpha) = -2 \lim_{\varepsilon \downarrow 0} \int_{-\infty}^\infty e^{-\varepsilon\lambda} \lambda^{-\frac{1}{2}} \sin(2\alpha\sqrt{\lambda}) \, d\rho(\lambda)$$

at any point, $\alpha$, of Lebesgue continuity for $q$. (1.23) is only proven for a restricted class of $q$'s including Case 1, 2 and those $q$'s satisfying

$$q(x) \geq -Cx^2, \qquad x \geq R$$

for some $R > 0$, $C > 0$, which are always in the limit point case at infinity. We will use (1.23) as our point of departure for relating $A(\alpha)$ to scattering data at the end of Section 8.

In order to prove (1.21) for finite $b$, we need to analyze the finite $b$ case extending (1.15) to all $a$ including $a = \infty$ (by allowing $A$ to have $\delta$ and $\delta'$ singularities at multiples of $b$). This was done in [35] for $\kappa$ real and positive and $a < \infty$. We now need results in the entire region $\text{Re}(\kappa) \geq K_0$, and this is what we do in Section 4. Explicitly, we will prove

THEOREM 1.6.  *In Case 1, there are $A_n, B_n$ for $n = 1, 2, \ldots,$ and a function $A(\alpha)$ on $(0, \infty)$ with*

(i) $|A_n| \leq C$.

(ii) $|B_n| \leq Cn$.

(iii) $\int_0^a |A(\alpha)| \, d\alpha \leq C \exp(K_0|a|)$ *so that for* $\text{Re}(\kappa) > \frac{1}{2}K_0$:

$$(1.24) \quad m(-\kappa^2) = -\kappa - \sum_{n=1}^\infty A_n \kappa e^{-2\kappa bn} - \sum_{n=1}^\infty B_n e^{-2\kappa bn} - \int_0^\infty A(\alpha) e^{-2\alpha\kappa} \, d\alpha.$$

In Section 6, we will use (1.21) to obtain *a priori* bounds on $\int_{-R}^0 d\rho(\lambda)$ as $R \to \infty$.

Section 9 includes further discussion of the significance of (1.21) and the connection between $A$ and the Gel′fand-Levitan transformation kernel.

Sections 10 and 11 present a few simple examples where one can compute $A$ explicitly. One of the examples, when combined with a general comparison theorem, allows us to prove the general bound
$$|A(\alpha)| \leq \alpha^{-1}\gamma(\alpha)e^{2\alpha\gamma(\alpha)},$$
where $\gamma(\alpha) = \sup_{0\leq x \leq \alpha} |q(x)|^{\frac{1}{2}}$ and this lets us extend (1.17) to bounded $q$.

In the appendix we discuss analogs of (1.15) for the other $m$-functions that arise in the Weyl-Titchmarsh theory.

While we will not discuss the theory in detail in this paper, we end this introduction by recalling the major thrust of [35] — the connection between $A$ and inverse theory (which holds for the principal $m$-function but not for the $m$-functions discussed in the appendix). Namely, there is an $A(\alpha, x)$ function associated to $m(z, x)$ by

$$(1.25) \qquad m(-\kappa^2, x) = -\kappa - \int_0^a A(\alpha, x) e^{-2\alpha\kappa} \, d\alpha + \tilde{O}(e^{-2a\kappa})$$

for $a < b - x$. This, of course, follows from Theorem 1.1 by translating the origin. The point is that $A$ satisfies the simple differential equation in the distributional sense

$$(1.26) \qquad \frac{\partial A}{\partial x}(\alpha, x) = \frac{\partial A}{\partial \alpha}(\alpha, x) + \int_0^\alpha A(\alpha - \beta, x) A(\beta, x) \, d\beta.$$

This is proven in [35] for $q \in L^1((0, a))$ (and some other $q$'s) and so holds in the generality of this paper since Theorem 1.3 implies $A(\alpha, x)$ for $\alpha + x \leq a$ is only a function of $q(y)$ for $y \in [0, a]$.

Moreover, by (1.16), we have

$$(1.27) \qquad \lim_{\alpha \downarrow 0} |A(\alpha, x) - q(\alpha + x)| = 0$$

uniformly in $x$ on compact subsets of the real line, so by the uniqueness theorem for solutions of (1.26) (proven in [35]), $A$ on $[0, a]$ determines $q$ on $[0, a]$.

In the limit circle case, there is an additional issue to discuss. Namely, that $m(z, x = 0)$ determines the boundary condition at $\infty$. This is because, as we just discussed, $m$ determines $A$ which determines $q$ on $[0, \infty)$. $m(z, x = 0)$ and $q$ determine $m(z, x)$ by the Riccati equation. Once we know $m$, we can recover $u(z = i, x) = \exp(\int_0^x m(z = i, y) \, dy)$, and so the particular solution that defined the boundary condition at $\infty$.

Thus, the inverse spectral theory aspects of the framework easily extend to the general case of potentials considered in the present paper.

With the exception of Theorem 2.1 for potentials $q \in L^1((0, \infty))$ of the first paper in this series [35], whose method of proof we follow in Section 4, we have made every effort to keep this paper independently readable and self-contained.

F.G. would like to thank C. Peck and T. Tombrello for the hospitality of Caltech where this work was done.

## 2. Using the Riccati equation

As explained in the introduction, the Riccati equation and *a priori* control on $m_j$ allow one to obtain exponentially small estimates on $m_1 - m_2$ (Theorem 1.5).

PROPOSITION 2.1. *Let $m_1(x), m_2(x)$ be two absolutely continuous functions on $[a, b]$ so that for some $Q \in L^1((a, b))$,*

(2.1) $$m_j'(x) = Q(x) - m_j(x)^2, \qquad j = 1, 2, \ x \in (a, b).$$

*Then*

$$[m_1(a) - m_2(a)] = [m_1(b) - m_2(b)] \exp\left(\int_a^b [m_1(y) + m_2(y)]\, dy\right).$$

*Proof.* Let $f(x) = m_1(x) - m_2(x)$ and $g(x) = m_1(x) + m_2(x)$. Then

$$f'(x) = -f(x)\, g(x),$$

from which it follows that

$$f(x) = f(b) \exp\left[\int_x^b g(y)\, dy\right]. \qquad \square$$

As an immediate corollary, we have the following (this implies Theorem 1.3)

THEOREM 2.2. *Let $m_j(x, -\kappa^2)$ be functions defined for $x \in [a, b]$ and $\kappa \in K$ some region of $\mathbb{C}$. Suppose that for each $\kappa$ in $K$, $m_j$ is absolutely continuous in $x$ and satisfies (N.B.: $q$ is the same for $m_1$ and $m_2$),*

$$m_j'(x, -\kappa^2) = q(x) + \kappa^2 - m_j(x, -\kappa^2)^2, \qquad j = 1, 2.$$

*Suppose $C$ is such that for each $x \in [a, b]$ and $\kappa \in K$,*

(2.2) $$|m_j(x, -\kappa^2) + \kappa| \leq C, \qquad j = 1, 2,$$

*then*

(2.3) $$|m_1(a, -\kappa^2) - m_2(a, -\kappa^2)| \leq 2C \exp[-2(b-a)[\operatorname{Re}(\kappa) - C]].$$

## 3. Atkinson's method

Theorem 2.2 places importance on *a priori* bounds of the form (2.2). Fortunately, by modifying ideas of Atkinson [3], we can obtain estimates of this form as long as $\operatorname{Im}(\kappa)$ is bounded away from zero.

Throughout this section, $b \leq \infty$ and $q \in L^1((0,a))$ for all $a < b$. For each $\kappa$ with $\text{Im}(\kappa) \neq 0$ and $\text{Re}(\kappa) > 0$, we suppose we are given a solution $u(x, -\kappa^2)$ of

(3.1) $$-u'' + qu = -\kappa^2 u,$$

which satisfies (note that $z = -\kappa^2$, so $\text{Im}(z) = -2\text{Re}(\kappa)\text{Im}(\kappa)$)

(3.2) $$-\text{Im}(\kappa)[\text{Im}(u'(x, -\kappa^2)/u(x, -\kappa^2))] > 0,$$

where $u' = \frac{\partial u}{\partial x}$. The examples to bear in mind are firstly $b < \infty$, $q \in L^1((0,b))$, and $u$ satisfies (3.1) with

$$u'(b_-, -\kappa^2) + hu(b_-, -\kappa^2) = 0 \qquad (|h| < \infty)$$

or

$$u(b_-, -\kappa^2) = 0 \qquad (h = \infty)$$

and secondly, $b = \infty$, and either $q$ limit point at infinity or $q$ limit circle with some boundary condition picked at $b$. Then take $u$ to be an $L^2$ solution of (3.1). In either case, $u$ can be chosen analytic in $\kappa$ although the bounds in Propositions 3.1 and 3.2 below do not require that.

Atkinson's method allows us to estimate $|m(-\kappa^2) + \kappa|$ in two steps. We will fix some $a < b$ finite and define $m_0(-\kappa^2)$ by solving

(3.3a) $$m_0'(-\kappa^2, x) = q(x) + \kappa^2 - m_0(-\kappa^2, x)^2,$$

(3.3b) $$m_0(-\kappa^2, a) = -\kappa$$

and then setting

(3.3c) $$m_0(-\kappa^2) := m_0(-\kappa^2, 0_+).$$

We will prove

PROPOSITION 3.1. *There is a $C > 0$ depending only on $q$ and a universal constant $E > 0$ so that if $\text{Re}(\kappa) \geq C$ and $\text{Im}(\kappa) \neq 0$, then*

(3.4) $$|m(-\kappa^2) - m_0(-\kappa^2)| \leq E \frac{|\kappa|^2}{|\text{Im}(\kappa)|} e^{-2a\text{Re}(\kappa)}.$$

*In fact, one can take*

$$C = \max\left(a^{-1}\ln(6), 4\int_0^a |q(x)|\,dx\right), \qquad E = \frac{3 \cdot 2 \cdot 12^2}{5}.$$

PROPOSITION 3.2. *There exist constants $D_1$ and $D_2$ (depending only on $a$ and $q$), so that for $\text{Re}(\kappa) > D_1$,*

$$|m_0(-\kappa^2) + \kappa| \leq D_2.$$

*Indeed, one can take*
$$D_1 = D_2 = 2\int_0^a |q(x)|\,dx.$$

These propositions together with Theorem 1.2 yield the following explicit form of Theorem 1.3.

THEOREM 3.3. *Let $q_1, q_2$ be defined on $(0, b_j)$ with $b_j > a$ for $j = 1, 2$. Suppose that $q_1 = q_2$ on $[0, a]$. Pick $\delta$ so that $a + \delta \leq \min(b_1, b_2)$ and let $\eta = \sup_{0 \leq x \leq a; j=1,2}(\int_x^{x+\delta} |q_j(y)|\,dy)$. Then if $\mathrm{Re}(\kappa) \geq \max(4\eta, \delta^{-1}\ln(6))$ and $\mathrm{Im}(\kappa) \neq 0$, we have that*
$$|m_1(-\kappa^2) - m_2(-\kappa^2)| \leq 2F(\kappa)\exp(-2a[\mathrm{Re}(\kappa) - F(\kappa)]),$$
*where*
$$F(\kappa) = 2\eta + \frac{864}{5}\frac{|\kappa|^2}{|\mathrm{Im}(\kappa)|}e^{-2\delta \mathrm{Re}(\kappa)}.$$

*Remarks.* 1. To obtain Theorem 1.3, we need only note that in the region $\arg(\kappa) \in (-\frac{\pi}{2} + \varepsilon, -\varepsilon)$, $|\kappa| \geq K_0$, $F(\kappa)$ is bounded.

2. We need not require that $\arg(\kappa) < -\varepsilon$ to obtain $F$ bounded. It suffices, for example, that $\mathrm{Re}(\kappa) \geq |\mathrm{Im}(\kappa)| \geq e^{-\alpha \mathrm{Re}(\kappa)}$ for some $\alpha < 2\delta$.

3. For $F$ to be bounded, we need not require that $\arg(\kappa) > -\frac{\pi}{2} + \varepsilon$. It suffices that $|\mathrm{Im}(\kappa)| \geq \mathrm{Re}(\kappa) \geq \alpha \ln[|\mathrm{Im}(\kappa)|]$ for some $\alpha > (2\delta)^{-1}$. Unfortunately, this does not include the region $\mathrm{Im}(-\kappa^2) = c$, $\mathrm{Re}(-\kappa^2) \to \infty$, where $\mathrm{Re}(\kappa)$ goes to zero as $|\kappa|^{-1}$. However, as $\mathrm{Re}(-\kappa^2) \to \infty$, we only need that $|\mathrm{Im}(-\kappa^2)| \geq 2\alpha|\kappa|\ln(|\kappa|)$.

As a preliminary to the proof of Proposition 3.1, we have

LEMMA 3.4. *Let $A, B, C, D \in \mathbb{C}$ so that $AD - BC = 1$ and so that $D \neq 0 \neq \mathrm{Im}(\frac{C}{D})$. Let $f$ be the fractional linear transformation*
$$f(\zeta) = \frac{A\zeta + B}{C\zeta + D}.$$
*Then $f[\mathbb{R} \cup \{\infty\}]$ is a circle of diameter*

(3.5) $$|D|^{-2}\left|\mathrm{Im}\left(\frac{C}{D}\right)\right|^{-1} = |\mathrm{Im}(\bar{C}D)|^{-1}.$$

*Remark.* If $|D| = 0$ or $\mathrm{Im}(\frac{C}{D}) = 0$, then $f[\mathbb{R} \cup \{\infty\}]$ is a straight line.

*Proof.* Consider first $g(\zeta) = \frac{\zeta}{a\zeta + 1} = \frac{1}{a + \zeta^{-1}}$. Then $g(0) = 0$ and $g'(0) = 1$, so $g[\mathbb{R} \cup \{\infty\}]$ is a circle tangent to the real axis. The other point on the imaginary axis has $\zeta = -\frac{1}{\mathrm{Re}(a)}$ with $g(-\frac{1}{\mathrm{Re}(a)}) = -\frac{i}{\mathrm{Im}(a)}$ so $\mathrm{diam}(g[\mathbb{R} \cup \{\infty\}]) = \frac{1}{|\mathrm{Im}(a)|}$.

Now write (using $AD - BC = 1$)

$$f(\zeta) = \frac{\zeta}{CD\zeta + D^2} + \frac{B}{D} = \frac{\zeta}{D^2[\frac{C}{D}\zeta + 1]} + \frac{B}{D}.$$

Thus letting $a = C/D$, $g(\zeta) = \frac{\zeta}{a\zeta+1}$ and writing $D = |D|e^{i\theta}$, we have that

$$f(\zeta) = e^{-2i\theta}|D|^{-2}g(\zeta) + \frac{B}{D}.$$

$B/D$ is a translation and $e^{-2i\theta}$ a rotation, and neither changes the diameter of a circle. So $\text{diam}(f[\mathbb{R} \cup \{\infty\}]) = |D|^{-2}\text{diam}(g[\mathbb{R} \cup \{\infty\}])$. □

Now let $\varphi(x, -\kappa^2), \theta(x, -\kappa^2)$ solve (3.1) with

(3.6a) $\qquad\qquad \varphi(0_+, -\kappa^2) = 0, \qquad \varphi'(0_+, -\kappa^2) = 1,$

(3.6b) $\qquad\qquad \theta(0_+, -\kappa^2) = 1, \qquad \theta'(0_+, -\kappa^2) = 0.$

Define

(3.7) $$f(\zeta) = -\frac{\theta(a, -\kappa^2)\zeta - \theta'(a, -\kappa^2)}{\varphi(a, -\kappa^2)\zeta - \varphi'(a, -\kappa^2)}.$$

LEMMA 3.5. *If $u$ solves (3.1) and $\frac{u'(a, -\kappa^2)}{u(a, -\kappa^2)} = \zeta$, then $\frac{u'(0, -\kappa^2)}{u(0, -\kappa^2)} = f(\zeta)$ with $f$ given by (3.7).*

*Proof.* Let $T = \begin{pmatrix} \varphi'(a, -\kappa^2) & \theta'(a, -\kappa^2) \\ \varphi(a, -\kappa^2) & \theta(a, -\kappa^2) \end{pmatrix}$. Then $T\begin{pmatrix} u'(0, -\kappa^2) \\ u(0, -\kappa^2) \end{pmatrix} = \begin{pmatrix} u'(a, -\kappa^2) \\ u(a, -\kappa^2) \end{pmatrix}$ by linearity of (3.1). By constancy of the Wronskian, $T$ has determinant 1 and thus

$$T^{-1} = \begin{pmatrix} \theta(a, -\kappa^2) & -\theta'(a, -\kappa^2) \\ -\varphi(a, -\kappa^2) & \varphi'(a, -\kappa^2) \end{pmatrix}$$

and so

$$\frac{u'(0, -\kappa^2)}{u(0, -\kappa^2)} = \frac{\theta(a, -\kappa^2)u'(a, -\kappa^2) - \theta'(a, -\kappa^2)u(a, -\kappa^2)}{-\varphi(a, -\kappa^2)u'(a, -\kappa^2) + \varphi'(a, -\kappa^2)u(a, -\kappa^2)}$$

$$= f\left(\frac{u'(a, -\kappa^2)}{u(a, -\kappa^2)}\right). \qquad □$$

COROLLARY 3.6.

(3.8) $$|m(-\kappa^2) - m_0(-\kappa^2)| \leq \frac{1}{|\text{Im}(\overline{\varphi(a, -\kappa^2)}\,\varphi'(a, -\kappa^2))|}.$$

*Proof.* We consider the case $\text{Im}(\kappa) < 0$. Let $m_1(-\kappa^2, x)$ be any solution of $m_1'(-\kappa^2, x) = q(x) + \kappa^2 - m_1^2(-\kappa^2, x)$. Then $\text{Im}(m_1'(-\kappa^2, x)) = 2\text{Re}(\kappa)\text{Im}(\kappa) - 2\text{Im}(m_1(-\kappa^2, x))\text{Re}(m_1(-\kappa^2, x))$. It follows that at a point where $\text{Im}(m_1) = 0$, that $\text{Im}(m_1') < 0$. Thus if $\text{Im}(m_1(-\kappa^2, y)) = 0$ for $y \in [0, a]$, then $\text{Im}(m_1(-\kappa^2, x)) < 0$ for $x \in (y, a]$. Thus $\text{Im}(m_1(-\kappa^2, a)) \geq 0$ implies $\text{Im}(m_1(-\kappa^2, 0)) > 0$, so $f$ maps $\mathbb{C}_+$ onto a circle in $\mathbb{C}_+$. Since $m_0(-\kappa^2, a) = -\kappa$ and $m(-\kappa^2, a)$ are in $\mathbb{C}_+$, both points are in $\mathbb{C}_+$ and so at $x = 0$, both lie inside the disc bounded by $f[\mathbb{R} \cup \{\infty\}]$. By $\det(T) = 1$ and Lemma 3.4, (3.8) holds. $\square$

*Proof of Proposition* 3.1. By (3.8), we need to estimate $\varphi(a, -\kappa^2)$. Define $w(x, -\kappa^2) = \text{Im}(\overline{\varphi(x, -\kappa^2)}\,\varphi'(x, -\kappa^2))$. Then, $w(0_+, -\kappa^2) = 0$ and by a standard Wronskian calculation, $w'(x, -\kappa^2) = -\text{Im}(-\kappa^2)\,|\varphi(x, -\kappa^2)|^2 = 2\text{Re}(\kappa)\text{Im}(\kappa)|\varphi(x, -\kappa^2)|^2$. Thus,

$$(3.9) \quad |\text{Im}(\overline{\varphi(a, -\kappa^2)}\,\varphi'(a, -\kappa^2))| = 2\text{Re}(\kappa)|\text{Im}(\kappa)|\int_0^a |\varphi(x, -\kappa^2)|^2\,dx.$$

$\varphi(x, -\kappa^2)$ satisfies the following integral equation ([5, §I.2]),

$$(3.10) \quad \varphi(x, -\kappa^2) = \frac{\sinh(\kappa x)}{\kappa} + \int_0^x \frac{\sinh(\kappa(x-y))}{\kappa}\,q(y)\varphi(y, -\kappa^2)\,dy.$$

Define $\beta(x, -\kappa^2) = \kappa e^{-x\text{Re}(\kappa)}\varphi(x, -\kappa^2)$. Then, by (3.10) and $|\sinh(\kappa\xi)| \leq e^{|\xi||\text{Re}(\kappa)|}$, $\xi \in \mathbb{R}$,

$$(3.11) \quad \left|\varphi(x, -\kappa^2) - \frac{\sinh(\kappa x)}{\kappa}\right| \leq \frac{e^{x\text{Re}(\kappa)}}{|\kappa|}\left[\sup_{0\leq y\leq x}|\beta(y, -\kappa^2)|\right]\frac{1}{|\kappa|}\int_0^x |q(y)|\,dy.$$

Moreover, (3.10) becomes

$$\beta(x, -\kappa^2) = e^{-x\text{Re}(\kappa)}\sinh(\kappa x) + \int_0^x \frac{\sinh(\kappa(x-y))}{\kappa}\,e^{-(x-y)\text{Re}(\kappa)}q(y)\beta(y, -\kappa^2)\,dy,$$

which implies that

$$(3.12) \quad |\beta(x, -\kappa^2)| \leq 1 + \int_0^x \frac{1}{|\kappa|}|q(y)|\,|\beta(y, -\kappa^2)|\,dy.$$

Pick $\kappa$ so that

$$(3.13) \quad |\kappa| \geq 4\int_0^a |q(y)|\,dy.$$

Then (3.12) implies

$$\sup_{0\leq x\leq a}|\beta(x, -\kappa^2)| \leq 1 + \frac{1}{4}\sup_{0\leq x\leq a}|\beta(x, -\kappa^2)|$$

so that

(3.14) $$\sup_{0 \leq x \leq a} |\beta(x, -\kappa^2)| \leq \frac{4}{3}.$$

Using (3.13) and (3.14) in (3.11), we get

(3.15) $$\left|\varphi(x, -\kappa^2) - \frac{\sinh(\kappa x)}{\kappa}\right| \leq \frac{1}{3} \frac{e^{x \operatorname{Re}(\kappa)}}{|\kappa|}.$$

Now $|\sinh(z)| \geq \sinh(|\operatorname{Re}(z)|) = \frac{1}{2}[e^{|\operatorname{Re}(z)|} - e^{-|\operatorname{Re}(z)|}]$, so (3.15) implies

(3.16) $$|\varphi(x, -\kappa^2)| \geq \frac{1}{6|\kappa|}[e^{x \operatorname{Re}(\kappa)} - 3e^{-x \operatorname{Re}(\kappa)}].$$

Now suppose

(3.17) $$a \operatorname{Re}(\kappa) \geq \ln(6).$$

Thus for $x \geq \frac{a}{2}$, (3.16) implies $|\varphi(x, -\kappa^2)| \geq \frac{1}{12|\kappa|} e^{x \operatorname{Re}(\kappa)}$ and we obtain

(3.18) $$\int_{\frac{a}{2}}^{a} |\varphi(y, -\kappa^2)|^2 \, dy \geq \frac{1}{288|\kappa|^2} \frac{1}{\operatorname{Re}(\kappa)} e^{2a \operatorname{Re}(\kappa)}[1 - e^{-a \operatorname{Re}(\kappa)}]$$
$$\geq \frac{5}{6} \frac{1}{288|\kappa|^2} \frac{1}{\operatorname{Re}(\kappa)} e^{2a \operatorname{Re}(\kappa)}.$$

Putting together (3.8), (3.9), and (3.18), we see that if (3.13) and (3.17) hold, then

$$|m(-\kappa^2) - m_0(-\kappa^2)| \leq \frac{3}{5} \times 288|\kappa|^2 \frac{1}{|\operatorname{Im}(\kappa)|} e^{-2a \operatorname{Re}(\kappa)}. \qquad \square$$

*Proof of Proposition* 3.2. Let

(3.19) $$\eta = \int_0^a |q(y)| \, dy.$$

Let $z(x, -\kappa^2)$ solve (3.1) with boundary conditions $z(a, -\kappa^2) = 1$, $z'(a, -\kappa^2) = -\kappa$, and let

$$\gamma(x, -\kappa^2) = \kappa + \frac{z'(x, -\kappa^2)}{z(x, -\kappa^2)}.$$

Then the Riccati equation for $m_0(-\kappa^2, a)$ becomes

(3.20) $$\gamma'(x, -\kappa^2) = q(x) - \gamma(x, -\kappa^2)^2 + 2\kappa \gamma(x, -\kappa^2)$$

and we have

(3.21) $$\gamma(a, -\kappa^2) = 0.$$

Thus $\gamma(x, -\kappa^2)$ satisfies

(3.22) $$\gamma(x, -\kappa^2) = -\int_x^a e^{-2\kappa(y-x)}[q(y) - \gamma^2(y, -\kappa^2)] \, dy.$$

Define $\Gamma(x, -\kappa^2) = \sup_{x \leq y \leq a} |\gamma(y, -\kappa^2)|$. Since $\text{Re}(\kappa) > 0$, (3.22) implies that

$$\Gamma(x, -\kappa^2) \leq \eta + \tfrac{1}{2}(\text{Re}(\kappa))^{-1}\Gamma(x, -\kappa^2)^2. \tag{3.23}$$

Suppose that

$$\text{Re}(\kappa) > 2\eta. \tag{3.24}$$

Then (3.23) implies

$$\Gamma(x, -\kappa^2) < \eta + \tfrac{1}{4}\eta^{-1}\Gamma(x, -\kappa^2)^2. \tag{3.25}$$

(3.25) implies $\Gamma(x, -\kappa^2) \neq 2\eta$. Since $\Gamma(a, -\kappa^2) = 0$ and $\Gamma$ is continuous, we conclude that $\Gamma(0_+, -\kappa^2) < 2\eta$, so $|\gamma(0_+, -\kappa^2)| < 2\eta$; hence $|(m_0(-\kappa^2) + \kappa| \leq 2\eta$. □

*Remark.* There is an interesting alternate proof of Proposition 3.2 that has better constants. It begins by noting that $m_0(-\kappa^2)$ is the $m$-function for the potential which is $q(x)$ for $x \leq a$ and $0$ for $x > a$. Thus Theorem 1.2 applies. So using the bounds (1.16) for $A$, we see immediately that for $\text{Re}(\kappa) > \tfrac{1}{2}\|q\|_1$,

$$|m_0(-\kappa^2) + \kappa| \leq \|q\|_1 + \frac{\|q\|_1^2}{[2\text{Re}(\kappa) - \|q\|_1]},$$

where $\|q\|_1 = \int_0^a |q(y)|\, dy$.

## 4. Finite $b$ representations with no errors

Theorem 1.2 says that if $b = \infty$ and $q \in L^1((0, \infty))$, then (1.17) holds, a Laplace transform representation for $m$ without errors. It is, of course, of direct interest that such a formula holds, but we are especially interested in a particular consequence of it — namely, that it implies that the formula (1.15) with error holds in the region $\text{Re}(\kappa) > K_0$ with error uniformly bounded in $\text{Im}(\kappa)$; that is, we are interested in

THEOREM 4.1.  *If $q \in L^1((0, \infty))$ and $\text{Re}(\kappa) > \tfrac{1}{2}\|q\|_1$, then for all $a$:*

$$\left| m(-\kappa^2) + \kappa + \int_0^a A(\alpha) e^{-2\alpha\kappa}\, d\alpha \right| \tag{4.1}$$
$$\leq \left[ \|q\|_1 + \frac{\|q\|_1^2 e^{a\|q\|_1}}{2\text{Re}(\kappa) - \|q\|_1} \right] e^{-2a\text{Re}(\kappa)}.$$

*Proof.* An immediate consequence of (1.17) and the estimate $|A(\alpha) - q(\alpha)| \leq \|q\|_1^2 \exp(\alpha \|q\|_1)$. □

Our principal goal in this section is to prove an analog of this result in the case $b < \infty$. To do so, we will need to first prove an analog of (1.17) in the case $b < \infty$ — something of interest in its own right. The idea will be to mimic the proof of Theorem 2 from [35] but use the finite $b$, $q^{(0)}(x) = 0$, $x \geq 0$ Green's function where [35] used the infinite $b$ Green's function. The basic idea is simple, but the arithmetic is a bit involved.

We will start with the $h = \infty$ case. Three functions for $q^{(0)}(x) = 0$, $x \geq 0$ are significant. First, the kernel of the resolvent $(-\frac{d^2}{dx^2} + \kappa^2)^{-1}$ with $u(0_+) = u(b_-) = 0$ boundary conditions. By an elementary calculation (see, e.g., [35, §5]), it has the form

$$(4.2) \qquad G^{(0)}_{h=\infty}(x, y, -\kappa^2) = \frac{\sinh(\kappa x_<)}{\kappa} \left[ \frac{e^{-\kappa x_>} - e^{-\kappa(2b - x_>)}}{1 - e^{-2\kappa b}} \right],$$

with $x_< = \min(x, y)$, $x_> = \max(x, y)$.

The second function is

$$(4.3) \qquad \psi^{(0)}_{h=\infty}(x, -\kappa^2) \equiv \lim_{y \downarrow 0} \frac{\partial G^{(0)}_{h=\infty}}{\partial y}(x, y, -\kappa^2) = \frac{e^{-\kappa x} - e^{-\kappa(2b - x)}}{1 - e^{-2\kappa b}}$$

and finally (notice that $\psi^{(0)}_{h=\infty}(0_+, -\kappa^2) = 1$ and $\psi^{(0)}_{h=\infty}$ satisfies the equations $-\psi'' = -\kappa^2 \psi$ and $\psi(b_-, -\kappa^2) = 0$):

$$(4.4) \qquad m^{(0)}_{h=\infty}(-\kappa^2) = \psi^{(0)\prime}_{h=\infty}(0_+, -\kappa^2) = -\frac{\kappa + \kappa e^{-2\kappa b}}{1 - e^{-2\kappa b}}.$$

In (4.4), prime means $d/dx$.

Fix now $q \in C_0^\infty((0, b))$. The pair of formulas

$$\left( -\frac{d^2}{dx^2} + q + \kappa^2 \right)^{-1} = \sum_{n=0}^\infty (-1)^n \left( -\frac{d^2}{dx^2} + \kappa^2 \right)^{-1} \left[ q \left( -\frac{d^2}{dx^2} + \kappa^2 \right)^{-1} \right]^n$$

and

$$m(-\kappa^2) = \lim_{x < y;\, y \downarrow 0} \frac{\partial^2 G(x, y, -\kappa^2)}{\partial x \partial y}$$

yields the following expansion for the $m$-function of $-\frac{d^2}{dx^2} + q$ with $u(b_-) = 0$ boundary conditions.

PROPOSITION 4.2. *Let $q \in C_0^\infty((0, b))$, $b < \infty$. Then*

$$(4.5) \qquad m(-\kappa^2) = \sum_{n=0}^\infty M_n(-\kappa^2; q),$$

*where*

$$(4.6) \qquad M_0(-\kappa^2; q) = m^{(0)}_{h=\infty}(-\kappa^2),$$

$$(4.7) \qquad M_1(-\kappa^2; q) = -\int_0^b q(x) \psi^{(0)}_{h=\infty}(x, -\kappa^2)^2\, dx,$$

and for $n \geq 2$,

$$
\begin{aligned}
(4.8)\quad M_n(-\kappa^2; q) = (-1)^n &\int_0^b dx_1 \ldots \int_0^b dx_n\, q(x_1)\ldots q(x_n) \\
&\times \psi_{h=\infty}^{(0)}(x_1, -\kappa^2)\psi_{h=\infty}^{(0)}(x_n, -\kappa^2) \prod_{j=1}^{n-1} G_{h=\infty}^{(0)}(x_j, x_{j+1}, -\kappa^2).
\end{aligned}
$$

The precise region of convergence is unimportant since we will eventually expand regions by analytic continuation. For now, we note it certainly converges in the region $\kappa$ real with $\kappa^2 > \|q\|_\infty$.

We want to write each term in (4.5) as a Laplace transform. We begin with (4.6), using (4.4)

$$
(4.9)\quad M_0(-\kappa^2; q) = -\kappa - \frac{2\kappa e^{-2\kappa b}}{1 - e^{-2\kappa b}} = -\kappa - 2\kappa \sum_{j=1}^\infty e^{-2j\kappa b}.
$$

Next, note by (4.3) that

$$
(4.10)\quad \psi_{h=\infty}^{(0)}(x, -\kappa^2) = \sum_{j=0}^\infty e^{-\kappa(x+2bj)} - \sum_{j=0}^\infty e^{-\kappa(2bj+(2b-x))},
$$

so

$$
(4.11)\quad \psi_{h=\infty}^{(0)}(x, -\kappa^2)^2 = \sum_{j=0}^\infty e^{-\kappa(2x+2bj)}(j+1) \\
+ \sum_{j=0}^\infty e^{-\kappa(2bj+(4b-2x))}(j+1) - 2\sum_{j=1}^\infty j e^{-2b\kappa j};
$$

hence,

$$
(4.12)\quad M_1(-\kappa^2; q) = 2\left[\int_0^b q(x)\,dx\right]\sum_{j=1}^\infty j e^{-2b\kappa j} - \int_0^\infty A_1(\alpha) e^{-2\alpha\kappa}\, d\alpha,
$$

where
$$
(4.13)\quad A_1(\alpha) = \begin{cases} q(\alpha), & 0 \leq \alpha < b, \\ (n+1)q(\alpha - nb) + nq((n+1)b - \alpha), & nb \leq \alpha < (n+1)b, \\ & n = 1, 2, \ldots. \end{cases}
$$

To manipulate $M_n$ for $n \geq 2$, we first rewrite (4.10) as

$$
(4.14)\quad \psi_{h=\infty}^{(0)}(x, -\kappa^2) = \sum_{j=0}^\infty \psi^{(0),(j)}(x, -\kappa^2)
$$

where

$$
(4.15)\quad \psi^{(0),(j)}(x, -\kappa^2) = (-1)^j \exp(-\kappa X_j(x)),
$$

with

(4.16) $$X_j(x) = \begin{cases} x + bj, & j = 0, 2, \ldots, \\ b - x + bj, & j = 1, 3, \ldots, \end{cases}$$

and then for $n \geq 2$

(4.17) $$M_n(-\kappa^2) = \sum_{j,p=0}^{\infty} M_{n,j,p}(-\kappa^2),$$

where

(4.18)
$$M_{n,j,p}(-\kappa^2) = (-1)^n \int_0^b dx_1 \ldots \int_0^b dx_n\, q(x_1) \ldots q(x_n)$$
$$\times \psi^{(0),(j)}(x_1, -\kappa^2) \psi^{(0),(p)}(x_n, -\kappa^2) \prod_{j=1}^{n-1} G_{h=\infty}^{(0)}(x_j, x_{j+1}, -\kappa^2).$$

Next use the representation from [35],

$$\frac{\sinh(\kappa x_<)}{\kappa} e^{-\kappa x_>} = \frac{1}{2} \int_{|x-y|}^{x+y} e^{-\kappa \ell}\, d\ell,$$

to rewrite (4.2) as

$$G_{h=\infty}^{(0)}(x, y, -\kappa^2) = \frac{1}{2} \int_{|x-y|}^{x+y} \left[ \frac{e^{-\kappa \ell} - e^{-\kappa(2b-\ell)}}{1 - e^{-2\kappa b}} \right] d\ell$$

$$= \frac{1}{2} \int_{S_+(x,y)} e^{-\kappa \ell}\, d\ell - \frac{1}{2} \int_{S_-(x,y)} e^{-\kappa \ell}\, d\ell,$$

where

$$S_+(x, y) = \bigcup_{n=0}^{\infty} [|x - y| + 2nb, x + y + 2nb]$$

and

$$S_-(x, y) = \bigcup_{n=0}^{\infty} [2b(n+1) - x - y, 2b(n+1) - |x - y|].$$

Each union consists of disjoint intervals although the two unions can overlap. The net result is that

(4.19) $$G_{h=\infty}^{(0)}(x, y, -\kappa^2) = \frac{1}{2} \int_0^{\infty} U(x, y, \ell) e^{-\kappa \ell}\, d\ell,$$

where $U$ is $+1$, $-1$, or $0$. The exact values of $U$ are complicated — that $|U| \leq 1$ is all we will need.

Plugging (4.19) in (4.18), we obtain

$$M_{n,j,p}(-\kappa^2) = \frac{(-1)^{n+j+p}}{2^{n-1}} \int_0^b dx_1 \ldots \int_0^b dx_n \int_0^\infty d\ell_1 \ldots \int_0^\infty d\ell_{n-1}$$

$$\times \; q(x_1)\ldots q(x_n) \prod_{j=1}^{n-1} U(x_j, x_{j+1}, \ell_j)$$

$$\times \; \exp(-\kappa[\ell_1 + \cdots + \ell_{n-1} + X_j(x_1) + X_p(x_n)]).$$

Letting $\alpha = \frac{1}{2}[\ell_1 + \cdots + \ell_{n-1} + X_j(x_1) + X_p(x_n)]$ and changing from $d\ell_{n-1}$ to $d\alpha$ (since $n \geq 2$, there is an $\ell_{n-2}$), we see that

$$(4.20) \qquad M_{n,j,p}(-\kappa^2) = -\int_{\frac{1}{2}b(j+p)}^\infty A_{n,j,p}(\alpha) e^{-2\alpha\kappa}\, d\alpha,$$

where
(4.21)
$$A_{n,j,p}(\alpha) = \frac{(-1)^{n+j+p}}{2^{n-2}} \int_0^b dx_1 \ldots \int_0^b dx_n \int_{R(x_1,\ldots x_n, \ell_1,\ldots,\ell_{n-2})} d\ell_1 \ldots d\ell_{n-2}$$

$$\times \prod_{j=1}^{n-2} U(x_j, x_{j+1}, \ell_j) U(x_{n-1}, x_n, 2\alpha - \ell_1 \ldots \ell_{n-2} - X_j(x_1) - X_p(x_n)),$$

where $R(x_1, \ldots, x_n, \ell_1, \ldots, \ell_{n-2})$ is the region

$$(4.22) \quad R(x_1, \ldots, x_n, \ell_1, \ldots, \ell_{n-2})$$
$$= \left\{ (\ell_1, \ldots, \ell_{n-2}) \;\bigg|\; \ell_i \geq 0 \text{ and } X_j(x_1) + X_p(x_n) + \sum_{k=1}^{n-2} \ell_k \leq 2\alpha \right\}.$$

In (4.20), the integral starts at $\frac{1}{2}b(j+p)$ since $\alpha \geq \frac{1}{2}[X_j(x_1) + X_p(x_n)]$ and (4.16) implies that $X_j(x) \geq bj$. For each value of $x$, $R$ is contained in the simplex $\{(\ell_1, \ldots, \ell_{n-2}) \mid \ell_i \geq 0 \text{ and } \sum_{k=1}^{n-2} \ell_k \leq 2\alpha\}$ which has volume $\frac{(2\alpha)^{n-2}}{(n-2)!}$. This fact and $|U| \leq 1$ employed in (4.21) imply

$$(4.23) \qquad |A_{n,j,p}(\alpha)| \leq \left(\int_0^b |q(x)|\, dx\right)^n \frac{\alpha^{n-2}}{(n-2)!}.$$

Moreover, by (4.20),

$$(4.24) \qquad A_{n,j,p}(\alpha) = 0 \quad \text{if } \alpha < \tfrac{1}{2} b(j+p).$$

For any fixed $\alpha$, the number of pairs $(j, p)$ with $j, p = 0, 1, 2 \ldots$ so that $\alpha > \frac{1}{2}b(j+p)$ is $\frac{1}{2}([\frac{2\alpha}{b}] + 1)([\frac{2\alpha}{b}] + 2)$; thus,

$$(4.25) \qquad M_n(-\kappa^2) = -\int_0^\infty A_n(\alpha) e^{-2\alpha\kappa}\, d\alpha,$$

with

$$(4.26) \qquad |A_n(\alpha)| \leq \left[\frac{(2\alpha + b)(2\alpha + 2b)}{2b^2}\right] \frac{\alpha^{n-2}}{(n-2)!} \|q\|_1^n.$$

As in [35], we can sum on $n$ from 2 to infinity and justify extending the result to all $q \in L^1((0,b))$. We therefore obtain

THEOREM 4.3 (Theorem 1.6 for $h = \infty$). *Let $b < \infty$, $h = \infty$, and $q \in L^1((0,b))$. Then for $\operatorname{Re}(\kappa) > \frac{1}{2}\|q\|_1$, we have that*

$$(4.27) \quad m(-\kappa^2) = -\kappa - \sum_{j=1}^{\infty} A_j \kappa e^{-2\kappa bj} - \sum_{j=1}^{\infty} B_j e^{-2\kappa bj} - \int_0^{\infty} A(\alpha) e^{-2\alpha\kappa}\, d\alpha,$$

*where*

(i) $A_j = 2$.

(ii) $B_j = -2j \int_0^b q(x)\, dx$.

(iii) $|A(\alpha) - A_1(\alpha)| \leq \frac{(2\alpha+b)(2\alpha+2b)}{2b^2} \|q\|_1^2 \exp(\alpha \|q\|_1)$ *with $A_1$ given by (4.13). In particular,*

$$\int_0^a |A(\alpha)|\, d\alpha \leq C(b, \|q\|_1)(1+a^2)\exp(a\|q\|_1).$$

As in the proof of Theorem 4.1, this implies

COROLLARY 4.4. *If $q \in L^1((0,\infty))$ and $\operatorname{Re}(\kappa) \geq \frac{1}{2}\|q\|_1 + \varepsilon$, then for all $a \in (0,b)$, $b < \infty$, we have that*

$$\left| m(-\kappa^2) + \kappa + \int_0^a A(\alpha) e^{-2\alpha\kappa}\, d\alpha \right| \leq C(a,\varepsilon) e^{-2a\operatorname{Re}(\kappa)},$$

*where $C(a,\varepsilon)$ depends only on $a$ and $\varepsilon$ (and $\|q\|_1$) but not on $\operatorname{Im}(\kappa)$.*

*Remark.* One can also prove results for $a > b$ if $b < \infty$ but this is the result we need in the next section.

The case $h = 0$ (Neumann boundary conditions at $b$) is almost the same. (4.2)–(4.4) are replaced by

$$(4.28) \quad G^{(0)}_{h=0}(x, y, -\kappa^2) = \frac{\sinh(\kappa x_<)}{\kappa}\left[\frac{e^{-\kappa x_>} + e^{-\kappa(2b-x_>)}}{1 + e^{-2\kappa b}}\right],$$

$$(4.29) \quad \psi^{(0)}_{h=0}(x, -\kappa^2) = \frac{e^{-\kappa x} + e^{-\kappa(2b-x)}}{1 + e^{-2\kappa b}},$$

$$(4.30) \quad m^{(0)}_{h=0}(-\kappa^2) = -\frac{\kappa - \kappa e^{-2\kappa b}}{1 + e^{-2\kappa b}}.$$

The only change in the further arguments is that $U$ can now take the values $0, \pm 1$, and $\pm 2$ so $|U| \leq 2$. That means that (4.26) becomes

$$|A_{n,h=0}(\alpha)| \leq 2\left[\frac{(2\alpha+b)(2\alpha+2b)}{2b^2}\right]\frac{(2\alpha)^{n-2}}{(n-2)!}\|q\|_1^n.$$

The net result is

THEOREM 4.5 (Theorem 1.6 for $h = 0$). *Let $b < \infty$, $h = 0$, and $q \in L^1((0,b))$. Then for $\mathrm{Re}(\kappa) > \|q\|_1$, (4.27) holds, where*

(i) $A_j = 2(-1)^j$.

(ii) $B_j = 2(-1)^{j+1} j \int_0^b q(x)\, dx$.

(iii) $|A(\alpha) - A_1(\alpha)| \leq \frac{(2\alpha+b)(2\alpha+2b)}{b^2} \|q\|_1^2 \exp(2\alpha\|q\|_1)$ *with $A_1$ given by*

$$A_{1,h=0}(\alpha) = \begin{cases} q(\alpha), & 0 \leq \alpha < b, \\ (-1)^n[(n+1)q(\alpha - nb) - nq((n+1)b - \alpha)], & nb \leq \alpha < (n+1)b, \\ & n = 1, 2, \ldots. \end{cases}$$

*In particular,*

$$\int_0^a |(A(\alpha)|\, d\alpha \leq C(b, \|q\|_1)(1 + a^2) \exp(2a\|q\|_1).$$

An analog of Corollary 4.4 holds, but we will wait for the general $h \in \mathbb{R}$ case to state it.

Finally, we turn to general $|h| < \infty$. In this case, (4.2)–(4.4) become

(4.31) $$G_h^{(0)}(x, y, -\kappa^2) = \frac{\sinh(\kappa x_<)}{\kappa} \psi_h^{(0)}(x_>, -\kappa^2),$$

(4.32) $$\psi_h^{(0)}(x, -\kappa^2) = \left[\frac{e^{-\kappa x} + \zeta(h, \kappa)e^{-\kappa(2b-x)}}{1 + \zeta(h, \kappa)e^{-2b\kappa}}\right],$$

(4.33) $$m_h^{(0)}(-\kappa^2) = -\kappa + 2\kappa \frac{\zeta(h, \kappa)e^{-2\kappa b}}{1 + \zeta(h, \kappa)e^{-2\kappa b}},$$

where

(4.34) $$\zeta(h, \kappa) = \frac{\kappa - h}{\kappa + h}.$$

To analyze this further, we need Laplace transform formulas for $\zeta$.

PROPOSITION 4.6. *The following formulas hold in the $\kappa$-region $h + \mathrm{Re}(\kappa) > 0$.*

(i) $\zeta(h, \kappa) = 1 - 4h \int_0^\infty e^{-\alpha(2\kappa+2h)}\, d\alpha$.

(ii) $\zeta(h, \kappa)^m = 1 + \sum_{j=1}^m (-1)^j \binom{m}{j} \frac{(4h)^j}{(j-1)!} \int_0^\infty \alpha^{j-1} e^{-\alpha(2\kappa+2h)}\, d\alpha$.

(iii) $\kappa\zeta(h, \kappa) = \kappa - 2h + 4h^2 \int_0^\infty e^{-\alpha(2\kappa+2h)}\, d\alpha$.

(iv) $\kappa\zeta(h, \kappa)^m = \kappa - 2mh - \frac{1}{4}\sum_{j=1}^m (-1)^j[\binom{m}{j} + 2\binom{m}{j+1}]\frac{(4h)^{j+1}}{(j-1)!} \int_0^\infty \alpha^{j-1} \times e^{-\alpha(2\kappa+2h)}\, d\alpha$, *where $\binom{m}{m+1}$ is interpreted as 0.*

*Proof.* Straightforward algebra. □

Rewriting (4.33) as

$$m_h^{(0)}(-\kappa^2) = -\kappa + 2\kappa \sum_{m=1}^{\infty} (-1)^{m+1} \zeta^m e^{-2m\kappa b}$$

and then using Proposition 4.6(iv), we find that

$$(4.35) \quad m_h^{(0)}(-\kappa^2) = -\kappa - 2\sum_{m=1}^{\infty}(-1)^m \kappa e^{-2m\kappa b}$$

$$- 4\sum_{m=1}^{\infty}(-1)^{m+1} mh e^{-2m\kappa b} - \int_{2b}^{\infty} A_{0,h}(\alpha) e^{-2\alpha\kappa}\, d\alpha,$$

where

$$(4.36) \quad A_{0,h}(\alpha) = \frac{1}{2}\sum_{m=1}^{\infty}(-1)^m \chi_{[2mb,\infty)}(\alpha) e^{-2(\alpha-2mb)h} \sum_{j=1}^{m}(-1)^j$$

$$\times \left[\binom{m}{j} + 2\binom{m}{j+1}\right](4h)^{j+1}\frac{(\alpha-2mb)^{j-1}}{(j-1)!}.$$

Using the crude estimates $(4h)^{j-1}(\alpha-2mb)^{j-1}\chi_{[2m,b)}(\alpha)/(j-1)! \leq \exp(4|h|\alpha)$, $\sum_{j=1}^{m}\binom{m}{j} \leq 2^m$, $\sum_{j=1}^{m}\binom{m}{j+1} \leq 2^m$, and $m \leq \alpha/2b$, we see that

$$(4.37) \quad |A_{0,h}(\alpha)| \leq \frac{3}{2}\left(\frac{\alpha}{2b}\right)\exp\left(\frac{\alpha}{2b}\ln(2)\right)\exp(6|h|\alpha).$$

A similar analysis of $\int_0^b q(x)\psi_{0,n}(x,-\kappa^2)^2\, dx$ shows that

(4.38)
$$-\int_0^b q(x)\psi_h^{(0)}(x,-\kappa^2)^2\, dx = -\left(\int_0^b q(x)\, dx\right) 2\sum_{m=1}^{\infty}(-1)^{m+1} m e^{-2b\kappa m}$$

$$- \int_0^b q(\alpha) e^{-2\alpha\kappa}\, d\alpha - \int_b^{\infty} A_{1,h}(\alpha) e^{-2\alpha\kappa}\, d\alpha,$$

where $A_{1,h}$ satisfies for suitable constants $C_1$ and $C_2$

(4.39)
$$|A_{1,h}(\alpha)| \leq C_1 \exp(C_2(|h|+1+b^{-1})\alpha)$$
$$\times [|q(\alpha-nb)| + |q((n+1)b - \alpha)|] \quad \text{for } nb \leq \alpha < (n+1)b.$$

Finally, using (4.31) and Proposition 4.6, we write

$$(4.40) \quad G_h^{(0)}(x,y,-\kappa^2) = \frac{1}{2}\int_0^{\infty} U(x,y,h,\ell) e^{-\kappa\ell}\, d\ell,$$

where

$$(4.41) \quad |U(x,y,h,\ell)| \leq C_3 \exp(C_4(|h|+1+b^{-1})\ell)$$

for suitable constants $C_3$ and $C_4$. From it, it follows that

$$\begin{aligned} M_n(-\kappa^2; q) &= (-1)^n \int_0^b dx_1 \ldots \int_0^b dx_n \, q(x_1) \ldots q(x_n) \\ &\quad \times \psi_h^{(0)}(x_1, -\kappa^2) \psi_h^{(0)}(x_n, -\kappa^2) \prod_{j=1}^{n-1} G_h^{(0)}(x_j, x_{j+1}, -\kappa^2) \\ &= -\int_0^\infty A_{n,h}(\alpha) e^{-2\kappa\alpha} \, d\alpha, \end{aligned}$$

where

$$|A_{n,h}(\alpha)| \leq C_5 \alpha^2 \exp(C_6(|h| + 1 + b^{-1})\alpha) \frac{\alpha^{n-2}}{(n-2)!} \|q\|_1^n, \quad n \geq 2.$$

We conclude

THEOREM 4.7 (Theorem 1.6 for general $|h| < \infty$). *Let $b < \infty$, $|h| < \infty$, and $q \in L^1((0, b))$. Then for $\mathrm{Re}(\kappa) > \frac{1}{2} D_1[\|q\|_1 + |h| + b^{-1} + 1]$ for a suitable universal constant $D_1$, (4.27) holds, where*

(i) $A_j = 2(-1)^j$.

(ii) $B_j = 2(-1)^{j+1} j[2h + \int_0^b q(x) \, dx]$.

(iii) $|A(\alpha) - q(\alpha)| \leq \|q\|_1^2 \exp(\alpha \|q\|_1)$ *if $|\alpha| < b$, and for any $a > 0$,*

$$\int_0^a |A(\alpha)| \, d\alpha \leq D_2(b, \|q\|_1, h) \exp(D_1 a(\|q\|_1 + |h| + b^{-1} + 1)).$$

Hence we immediately get

COROLLARY 4.8. *Fix $b < \infty$, $q \in L^1((0, b))$, and $|h| < \infty$. Fix $a < b$. Then there exist positive constants $C$ and $K_0$ so that for all complex $\kappa$ with $\mathrm{Re}(\kappa) > K_0$,*

$$\left| m(-\kappa^2) + \kappa + \int_0^a A(\alpha) e^{-2\alpha\kappa} \, d\alpha \right| \leq C e^{-2a\kappa}.$$

## 5. The relation between $A$ and $\rho$: Distributional form, I.

Our primary goal in the next five sections is to discuss a formula which formally says that

(5.1) $$A(\alpha) = -2 \int_{-\infty}^\infty \lambda^{-\frac{1}{2}} \sin(2\alpha\sqrt{\lambda}) \, d\rho(\lambda),$$

where for $\lambda \leq 0$, we define

$$\lambda^{-\frac{1}{2}} \sin(2\alpha\sqrt{\lambda}) = \begin{cases} 2\alpha & \text{if } \lambda = 0, \\ (-\lambda)^{-\frac{1}{2}} \sinh(2\alpha\sqrt{-\lambda}) & \text{if } \lambda < 0. \end{cases}$$

In a certain sense which will become clear, the left-hand side of (5.1) should be $A(\alpha) - A(-\alpha) + \delta'(\alpha)$.

To understand (5.1) at a formal level, note the basic formulas,

$$(5.2) \qquad m(-\kappa^2) = -\kappa - \int_0^\infty A(\alpha)e^{-2\alpha\kappa}\,d\alpha,$$

$$(5.3) \qquad m(-\kappa^2) = \operatorname{Re}(m(i)) + \int_{-\infty}^\infty \left[\frac{1}{\lambda+\kappa^2} - \frac{\lambda}{1+\lambda^2}\right] d\rho(\lambda),$$

and

$$(5.4) \qquad (\lambda+\kappa^2)^{-1} = 2\int_0^\infty \lambda^{-\frac{1}{2}} \sin(2\alpha\sqrt{\lambda})e^{-2\alpha\kappa}\,d\alpha,$$

which is an elementary integral if $\kappa > 0$ and $\lambda > 0$. Plug (5.4) into (5.3), formally interchange order of integrations, and (5.2) should only hold if (5.1) does. However, a closer examination of this procedure reveals that the interchange of order of integrations is not justified and indeed (5.1) is not true as a simple integral since, as we will see in the next section, $\int_0^R d\rho(\lambda) \underset{R\to\infty}{\sim} \frac{2}{3\pi}R^{\frac{3}{2}}$, which implies that (5.1) is not absolutely convergent. We will even see (in §9) that the integral sometimes fails to be conditionally convergent.

Our primary method for understanding (5.1) is as a distributional statement, that is, it will hold when smeared in $\alpha$ for $\alpha$ in $(0, b)$. We prove this in this section if $q \in L^1((0,\infty))$ or if $b < \infty$. In Section 7, we will extend this to all $q$ (i.e., all Cases 1–4) by a limiting argument using estimates we prove in Section 6. The estimates themselves will come from (5.1)! In Section 8, we will prove (5.1) as a pointwise statement where the integral is defined as an Abelian limit. Again, estimates from Section 6 will play a role.

Suppose $b < \infty$ or $b = \infty$ and $q \in L^1((0,b))$. Fix $a < b$ and $f \in C_0^\infty((0,a))$. Define

$$(5.5) \qquad m_a(-\kappa^2) := -\kappa - \int_0^a A(\alpha)e^{-2\alpha\kappa}\,d\alpha$$

for $\operatorname{Re}(\kappa) \geq 0$. Fix $\kappa_0$ real and let

$$g(y, \kappa_0, a) := m_a(-(\kappa_0 + iy)^2),$$

with $\kappa_0, a$ as real parameters and $y \in \mathbb{R}$ a variable. As usual, define the Fourier transform by (initially for smooth functions and then by duality for tempered distributions [33, Ch. IX])

$$(5.6) \qquad \hat{F}(k) = \frac{1}{\sqrt{2\pi}} \int_\mathbb{R} e^{-iky} F(y)\,dy, \qquad \check{F}(k) = \frac{1}{\sqrt{2\pi}} \int_\mathbb{R} e^{iky} F(y)\,dy.$$

Then by (5.5),

$$(5.7) \qquad \widehat{g}(k, \kappa_0, a) = -\sqrt{2\pi}\,\kappa_0 \delta(k) - \sqrt{2\pi}\,\delta'(k) - \frac{\sqrt{2\pi}}{2} e^{-k\kappa_0} A\left(\frac{k}{2}\right)\chi_{(0,2a)}(k).$$

Thus, since $f(0_+) = f'(0_+) = 0$, in fact, $f$ has support away from 0 and $a$,

$$
(5.8) \quad \int_0^a A(\alpha) f(\alpha) \, d\alpha = -\frac{2}{\sqrt{2\pi}} \int_0^a \overline{\tilde{g}}(2\alpha, \kappa_0, a) e^{2\alpha \kappa_0} f(\alpha) \, d\alpha
$$

$$
= -\frac{1}{\sqrt{2\pi}} \int_0^{2a} \overline{\tilde{g}}(\alpha, \kappa_0, a) e^{\alpha \kappa_0} f\left(\frac{\alpha}{2}\right) d\alpha
$$

$$
= -\frac{1}{\sqrt{2\pi}} \int_{\mathbb{R}} g(y, \kappa_0, a) \check{F}(y, \kappa_0) \, dy,
$$

where we have used the unitarity of $\hat{\ }$ and

$$
(5.9) \quad \check{F}(y, \kappa_0) = \frac{1}{\sqrt{2\pi}} \int_0^{2a} e^{\alpha(\kappa_0 + iy)} f\left(\frac{\alpha}{2}\right) d\alpha
$$

$$
= \frac{2}{\sqrt{2\pi}} \int_0^a e^{2\alpha(\kappa_0 + iy)} f(\alpha) \, d\alpha.
$$

Notice that

$$
(5.10) \quad |\check{F}(y, \kappa_0)| \leq C e^{2(a-\varepsilon)\kappa_0} (1 + |y|^2)^{-1}
$$

since $f$ is smooth and supported in $(0, a - \varepsilon)$ for some $\varepsilon > 0$.

By Theorem 4.1 and Corollary 4.8,

$$
(5.11) \quad |m_a(-(\kappa_0 + iy)^2) - m(-(\kappa_0 + iy)^2)| \leq C e^{-2a\kappa_0}
$$

for large $\kappa_0$, uniformly in $y$. From (5.8), (5.10), and (5.11), one concludes that

LEMMA 5.1. *Let $f \in C_0^\infty((0, a))$ with $0 < a < b$ and $q \in L^1((0, b))$. Then*

$$
(5.12)
$$
$$
\int_0^a A(\alpha) f(\alpha) \, d\alpha = \lim_{\kappa_0 \uparrow \infty} \left[ -\frac{1}{\pi} \int_{\mathbb{R}} m(-(\kappa_0 + iy)^2) \left[ \int_0^a e^{2\alpha(\kappa_0 + iy)} f(\alpha) \, d\alpha \right] dy \right].
$$

As a function of $y$, for $\kappa_0$ fixed, the alpha integral is $O((1 + y^2)^{-N})$ for all $N$ because $f$ is $C^\infty$. Now define

$$
(5.13) \quad \tilde{m}_R(-\kappa^2) = \left[ c_R + \int_{\lambda \leq R} \frac{d\rho(\lambda)}{\lambda + \kappa^2} \right],
$$

where $c_R$ is chosen so that $\tilde{m}_R \xrightarrow[R \to \infty]{} m$. Because $\int_{\mathbb{R}} \frac{d\rho(\lambda)}{1+\lambda^2} < \infty$, the convergence is uniform in $y$ for $\kappa_0$ fixed and sufficiently large. Thus in (5.12) we can replace $m$ by $m_R$ and take a limit (first $R \to \infty$ and then $\kappa_0 \uparrow \infty$). Since $f(0_+) = 0$, the $\int dy \, c_R \, d\alpha$-integrand is zero. Moreover, we can now interchange the $dy \, d\alpha$ and $d\rho(\lambda)$ integrals. The result is that

(5.14)
$$\int_0^a A(\alpha)f(\alpha)\,d\alpha = \lim_{\kappa_0\uparrow\infty}\lim_{R\to\infty}\int_{\lambda\leq R} d\rho(\lambda)$$
$$\times\left[\int_0^a d\alpha\, e^{2\alpha\kappa_0} f(\alpha)\left[-\frac{1}{\pi}\int_{\mathbb{R}}\frac{e^{2\alpha iy}\,dy}{(\kappa_0+iy)^2+\lambda}\right]\right].$$

In the case at hand, $d\rho$ is bounded below, say $\lambda \geq -K_0$. As long as we take $\kappa_0 > K_0$, the poles of $(\kappa_0+iy)^2+\lambda$ occur in the upper half-plane
$$y_\pm = i\kappa_0 \pm \sqrt{\lambda}.$$

Closing the contour in the upper plane, we find that if $\lambda \geq -K_0$,
$$-\frac{1}{\pi}\int_{\mathbb{R}}\frac{e^{2\alpha iy}\,dy}{(\kappa_0+iy)^2+\lambda} = -2e^{-2\alpha\kappa_0}\frac{\sin(2\alpha\sqrt{\lambda})}{\sqrt{\lambda}}.$$

Thus (5.14) becomes
$$\int_0^a A(\alpha)f(\alpha)\,d\alpha = -2\lim_{\kappa_0\uparrow\infty}\lim_{R\to\infty}\int_{\lambda\leq R}\left[\int_0^a f(\alpha)\frac{\sin(2\alpha\sqrt{\lambda})}{\sqrt{\lambda}}\,d\alpha\right]d\rho(\lambda).$$

$\kappa_0$ has dropped out and the $\alpha$ integral is bounded by $C(1+\lambda^2)^{-1}$, so we take the limit as $R \to \infty$ since $\int_{\mathbb{R}}\frac{d\rho(\lambda)}{1+\lambda^2} < \infty$. We have therefore proven the following result.

THEOREM 5.2. *Let $f \in C_0^\infty((0,a))$ with $a < b$ and either $b < \infty$ or $q \in L^1((0,\infty))$ with $b = \infty$. Then*

(5.15)
$$\int_0^a A(\alpha)f(\alpha)\,d\alpha = -2\int_{\mathbb{R}}\left[\int_0^a f(\alpha)\frac{\sin(2\alpha\sqrt{\lambda})}{\sqrt{\lambda}}\,d\alpha\right]d\rho(\lambda).$$

We will need to strengthen this in two ways. First, we want to allow $a > b$ if $b < \infty$. As long as $A$ is interpreted as a distribution with $\delta$ and $\delta'$ functions at $\alpha = nb$, this is easy. We also want to allow $f$ to have a nonzero derivative at $\alpha = 0$. The net result is

THEOREM 5.3. *Let $f \in C_0^\infty(\mathbb{R})$ with $f(-\alpha) = -f(\alpha)$, $\alpha \in \mathbb{R}$ and either $b < \infty$ or $q \in L^1((0,\infty))$ with $b = \infty$. Then*

(5.16)
$$-2\int_{\mathbb{R}}\left[\int_{-\infty}^\infty f(\alpha)\frac{\sin(2\alpha\sqrt{\lambda})}{\sqrt{\lambda}}\,d\alpha\right]d\rho(\lambda) = \int_{-\infty}^\infty \tilde{A}(\alpha)f(\alpha)\,d\alpha,$$

*where $\tilde{A}$ is the distribution*

(5.17a) $\tilde{A}(\alpha) = \chi_{(0,\infty)}(\alpha)A(\alpha) - \chi_{(-\infty,0)}(\alpha)A(-\alpha) + \delta'(\alpha)$

*if $b = \infty$ and*

$$(5.17\text{b}) \qquad \tilde{A}(\alpha) = \chi_{(0,\infty)}(\alpha)A(\alpha) - \chi_{(-\infty,0)}(\alpha)A(-\alpha) + \delta'(\alpha)$$

$$+ \sum_{j=1}^{\infty} B_j[\delta(\alpha - 2bj) - \delta(\alpha + 2bj)]$$

$$+ \sum_{j=1}^{\infty} \tfrac{1}{2} A_j[\delta'(\alpha - 2bj) + \delta'(\alpha + 2bj)]$$

if $b < \infty$, where $A_j, B_j$ are $h$ dependent and given in Theorems 4.3, 4.5, and 4.7.

The proof is identical to the argument above. $f(0)$ is still 0 but since $f'(0) \neq 0$, we carry it along.

*Example.* Let $b = \infty$, $q^{(0)}(x) = 0$, $x \geq 0$. Then $d\rho^{(0)}(\lambda) = \tfrac{1}{\pi}\chi_{[0,\infty)}(\lambda) \times \sqrt{\lambda}\, d\lambda$. Thus,

$$(5.18) \quad -2\int_{-\infty}^{\infty}\left[\int_{-\infty}^{\infty} f(\alpha)\frac{\sin(2\alpha\sqrt{\lambda})}{\sqrt{\lambda}}\,d\alpha\right]d\rho^{(0)}(\lambda)$$

$$= -\frac{2}{\pi}\int_0^{\infty}\left[\int_{-\infty}^{\infty} f(\alpha)\sin(2\alpha\sqrt{\lambda})\,d\alpha\right]d\lambda.$$

Next, change variables by $k = 2\sqrt{\lambda}$, that is, $\lambda = \tfrac{k^2}{4}$, and then change from $\int_0^{\infty} dk$ to $\tfrac{1}{2}\int_{-\infty}^{\infty} dk$ to obtain (recall $f(-\alpha) = -f(\alpha)$)

$$(5.18) = -\frac{1}{2\pi}\int_{-\infty}^{\infty}\left[\int_{-\infty}^{\infty} f(\alpha)\sin(\alpha k)\,d\alpha\right] k\,dk$$

$$= \frac{1}{\sqrt{2\pi}}\int_{-\infty}^{\infty} ik\,\check{f}(k)\,dk$$

$$= -f'(0)$$

$$= \int_{-\infty}^{\infty} f(\alpha)\delta'(\alpha)\,d\alpha,$$

as claimed in (5.16) and (5.17a) since $A^{(0)}(\alpha) = 0$, $\alpha \geq 0$.

## 6. Bounds on $\int_0^{\pm R} d\rho(\lambda)$

As we will see, (1.13) implies asymptotic results on $\int_{-R}^{R} d\rho(\lambda)$ and (5.1) will show that $\int_{-\infty}^{0} e^{b\sqrt{-\lambda}} d\rho(\lambda) < \infty$ for all $b > 0$ and more (for remarks on the history of the subject, see the end of this section). It follows from (5.3) that

$$\text{Im}(m(ia)) = a\int_{\mathbb{R}} \frac{d\rho(\lambda)}{\lambda^2 + a^2}, \qquad a > 0.$$

Thus, Everitt's result (1.13) (which also follows from our results in §§2 and 3) implies that
$$\lim_{a\to\infty} a^{\frac{1}{2}} \int_{\mathbb{R}} \frac{d\rho(\lambda)}{\lambda^2 + a^2} = 2^{-\frac{1}{2}}.$$

Standard Tauberian arguments (see, e.g., [34, §III.10], which in this case shows that on even functions $R^{\frac{3}{2}} d\rho(\frac{\lambda}{R}) \underset{R\to\infty}{\to} \frac{1}{2}\pi^{-1} |\lambda|^{\frac{1}{2}} d\lambda)$ then imply

THEOREM 6.1.

(6.1) $$\lim_{R\to\infty} R^{-\frac{3}{2}} \int_{-R}^{R} d\rho(\lambda) = \frac{2}{3\pi}.$$

*Remarks.* 1. This holds in all cases (1–4) we consider here, including some with supp$(d\rho)$ unbounded below.

2. Since we will see $\int_{-\infty}^{0} d\rho$ is bounded, we can replace $\int_{-R}^{R}$ by $\int_{0}^{R}$ in (6.1).

We will need the following *a priori* bound that follows from Propositions 3.1 and 3.2

PROPOSITION 6.2. *Let $d\rho$ be the spectral measure for a Schrödinger operator in Cases 1–4. Fix $a < b$. Then there is a constant $C_a$ depending only on $a$ and $\int_{0}^{a} |q(y)| dy$ so that*

(6.2) $$\int_{\mathbb{R}} \frac{d\rho(\lambda)}{1+\lambda^2} \leq C_a.$$

*Proof.* By Propositions 3.1 and 3.2, we can find $C_1$ and $z_1 \in \mathbb{C}_+$ depending only on $a$ and $\int_{0}^{a} |q(y)| dy$ so that
$$|m(z_1)| \leq C_1.$$
Thus,
$$\int_{\mathbb{R}} \frac{d\rho(\lambda)}{(\lambda - \text{Re}(z_1))^2 + (\text{Im}(z_1))^2} = \frac{\text{Im}(m(z_1))}{\text{Im}(z_1)} \leq \frac{C_1}{\text{Im}(z_1)},$$
and
$$\int_{\mathbb{R}} \frac{d\rho}{1+\lambda^2} \leq \frac{C_1}{\text{Im}(z_1)} \sup_{\lambda \in \mathbb{R}} \left[ \frac{(\lambda - \text{Re}(z_1))^2 + (\text{Im}(z_1))^2}{1+\lambda^2} \right] \equiv C_a. \quad \square$$

Our main goal in the rest of this section will be to bound $\int_{-\infty}^{0} e^{2\alpha\sqrt{-\lambda}} \times d\rho(\lambda)$ for any $\alpha < b$ and to find an explicit bound in terms of $\sup_{0 \leq x \leq \alpha+1} [-q(y)]$ when that sup is finite. As a preliminary, we need the following result from the standard limit circle theory [6, §9.4].

PROPOSITION 6.3. *Let $b = \infty$ and let $d\rho$ be the spectral measure for a problem of types 2–4. Let $d\rho_{R,h}$ be the spectral measure for the problem with $b = R < \infty, h$, and potential equal to $q(x)$ for $x \leq R$. Then there exists $h(R)$ so that*

$$d\rho_{R,h(R)} \underset{R \to \infty}{\to} d\rho,$$

*when smeared with any function $f$ of compact support. In particular, if $f \geq 0$, then*

$$\int_{\mathbb{R}} f(\lambda)\, d\rho(\lambda) \leq \overline{\lim_{R \to \infty}} \int_{\mathbb{R}} f(\lambda)\, d\rho_{R,h(R)}(\lambda).$$

This result implies that we need only obtain bounds for $b < \infty$ (where we have already proven (5.15)).

LEMMA 6.4. *If $\rho_1$ has support in $[-E_0, \infty)$, $E_0 > 0$, then*

$$\tag{6.3} \int_{-\infty}^{0} e^{\gamma\sqrt{-\lambda}}\, d\rho_1(\lambda) \leq e^{\gamma\sqrt{E_0}}(1 + E_0^2) \int_{-\infty}^{0} \frac{d\rho_1(\lambda)}{1 + \lambda^2}.$$

*Proof.* Obvious. □

Now let $f$ be fixed in $C_0^\infty((0,1))$ with $f \geq 0$ and $\int_0^1 f(y)\, dy = 1$. Let $f_{\alpha_0}(\alpha) = f(\alpha - \alpha_0)$. Let $d\rho_2$ be the spectral measure for some problem with $b \geq \alpha_0 + 1$ and let $d\rho_1$ be the spectral measure for the problem with $b = \alpha_0 + 1$, $h = \infty$, and the same potential on $[0, \alpha_0 + 1]$. Then, by Theorem 1.3, $A_1(\alpha) = A_2(\alpha)$ for $\alpha \in [0, \alpha_0 + 1]$ so $\int_{\alpha_0}^{\alpha_0+1} f_{\alpha_0}(\alpha)[A_1(\alpha) - A_2(\alpha)]\, d\alpha = 0$, and thus by Theorem 5.2,

$$\tag{6.4} \int_{\mathbb{R}} G_{\alpha_0}(\lambda)[d\rho_1(\lambda) - d\rho_2(\lambda)] = 0,$$

where

$$\tag{6.5} G_{\alpha_0}(\lambda) = \int_{\alpha_0}^{\alpha_0+1} f_{\alpha_0}(\alpha) \frac{\sin(2\alpha\sqrt{\lambda})}{\sqrt{\lambda}}\, d\alpha.$$

LEMMA 6.5. (i) *For $\lambda \geq 0$, $|G_{\alpha_0}(\lambda)| \leq 2(1 + \alpha_0)$.*

(ii) $|G_{\alpha_0}(\lambda)| \leq \lambda^{-2} \frac{1}{8} \int_0^1 |f'''(u)|\, du := C_0 \lambda^{-2}$ *for $\lambda > 0$.*

(iii) *For $\lambda \leq 0$, $|G_{\alpha_0}(\lambda)| \leq 2(\alpha_0 + 1)e^{2(\alpha_0+1)\sqrt{-\lambda}}$.*

(iv) *For $\lambda \leq 0$, $G_{\alpha_0}(\lambda) \geq \frac{1}{2\sqrt{-\lambda}}[e^{2\alpha_0\sqrt{-\lambda}} - 1]$.*

*Proof.* (i) Since $|\sin(x)| \leq |x|$, $|\sin(2\alpha\sqrt{\lambda})/\sqrt{\lambda}| \leq 2\alpha$. Thus, since $\text{supp}(f_{\alpha_0}) \subset [\alpha_0, \alpha_0 + 1]$ and $\int_{\alpha_0}^{\alpha_0+1} f_{\alpha_0}(\alpha)\, d\alpha = 1$, $|G_{\alpha_0}(\lambda)| \leq 2(1 + \alpha_0)$.

(ii) $\frac{1}{(2\lambda^{\frac{1}{2}})^3} \frac{d^3}{d\alpha^3} \cos(2\alpha\sqrt{\lambda}) = \sin(2\alpha\sqrt{\lambda})$, so this follows upon integrating by parts repeatedly.

(iii), (iv) For $y \geq 0$,
$$\frac{\sinh(y)}{y} = \frac{1}{y}\int_0^y \cosh(u)\,du$$
so $\frac{1}{2}e^u \leq \cosh u \leq e^u \leq e^y$, $0 \leq u \leq y$ implies
$$\frac{e^y - 1}{2y} \leq \frac{\sinh(y)}{y} \leq e^y.$$
This implies (iii) and (iv) given $\mathrm{supp}(f_{\alpha_0}) \subset [\alpha_0, \alpha_0 + 1]$, $f_{\alpha_0}(\alpha) \geq 0$, and $\int_{\alpha_0}^{\alpha_0+1} f_{\alpha_0}(\alpha)\,d\alpha = 1$. □

We can plug in these estimates into (6.4) to obtain
$$\int_{-\infty}^0 \frac{1}{2\sqrt{-\lambda}}\left[e^{2\alpha_0\sqrt{-\lambda}} - 1\right] d\rho_2(\lambda) \leq T_1 + T_2 + T_3,$$
where,
$$T_j = \max(4(1+\alpha_0), 2C_0) \int_{\mathbb{R}} \frac{d\rho_j(\lambda)}{1+\lambda^2}, \qquad j = 1, 2,$$
$$T_3 = \int_{-\infty}^0 2(\alpha_0 + 1)e^{2(\alpha_0+1)\sqrt{-\lambda}}\,d\rho_1(\lambda),$$
and we have used
$$1 \leq \frac{2}{1+\lambda^2}, \qquad 0 \leq \lambda \leq 1,$$
$$\frac{1}{\lambda^2} \leq \frac{2}{1+\lambda^2}, \qquad \lambda \geq 1.$$

Thus, Propositions 6.2, 6.3 and Lemma 6.4 together with
$$\frac{e^u - 1}{u} = \int_0^1 e^{yu}\,dy \geq e^{(1-\delta)u}\delta$$
for any $u > 0$ and any $\delta \in \mathbb{R}$ imply

THEOREM 6.6.  *Let $\rho$ be the spectral measure for some problem of the types 2–4. Let*

(6.6)
$$E(\alpha_0) := -\inf\left\{\int_0^{\alpha_0+1}(|\varphi_n'(x)|^2 + q(x)|\varphi(x)|^2\,dx)\,\bigg|\,\varphi \in C_0^\infty((0, \alpha_0+1)),\right.$$
$$\left.\int_0^{\alpha_0+1} |\varphi(x)|^2\,dx \leq 1\right\}.$$

*Then for all $\delta > 0$ and $\alpha_0 > 0$,*
(6.7)
$$\alpha_0\delta \int_{-\infty}^0 e^{2(1-\delta)\alpha_0\sqrt{-\lambda}}\,d\rho(\lambda) \leq \left[C_1(1+\alpha_0) + C_2(1+E(\alpha_0)^2)e^{2(\alpha_0+1)\sqrt{E(\alpha_0)}}\right],$$

where $C_1, C_2$ only depend on $\int_0^1 |q(x)|\, dx$. In particular,

$$\text{(6.8)} \qquad \int_{-\infty}^{0} e^{B\sqrt{-\lambda}}\, d\rho(\lambda) < \infty$$

for all $B < \infty$.

As a special case, suppose $q(x) \geq -C(x+1)^2$. Then $E(\alpha_0) \geq -C(\alpha_0+2)^2$ and we see that

$$\text{(6.9)} \qquad \int_{-\infty}^{0} e^{B\sqrt{-\lambda}} d\rho(\lambda) \leq D_1 e^{D_2 B^2}.$$

This implies

THEOREM 6.7. *If $d\rho$ is the spectral measure for a potential which satisfies*

$$\text{(6.10)} \qquad q(x) \geq -Cx^2, \qquad x \geq R$$

*for some $R > 0$, $C > 0$, then for $\varepsilon > 0$ sufficiently small,*

$$\text{(6.11)} \qquad \int_{-\infty}^{0} e^{-\varepsilon \lambda}\, d\rho(\lambda) < \infty.$$

*Remarks.* 1. Our proof shows in terms of the $D_2$ of (6.9), one only needs that $\varepsilon < 1/4D_2$.

2. Our proof implies that if

$$\lim_{x \to \infty} \frac{1}{x^2} \max(0, -q(x)) = 0,$$

then (6.11) holds for all $\varepsilon > 0$.

*Proof.* (6.9) implies that

$$\int_{-(n+1)^2}^{-n^2} d\rho(\lambda) \leq D_1 e^{D_2 B^2} e^{-Bn}.$$

Taking $B = n/2D_2$, we see that

$$\int_{-(n+1)^2}^{-n^2} d\rho(\lambda) \leq D_1 e^{-\frac{n^2}{4D_2}}.$$

Thus,

$$\int_{-\infty}^{0} e^{-\varepsilon \lambda}\, d\rho(\lambda) \leq \sum_{n=0}^{\infty} e^{\varepsilon(n+1)^2} \int_{-(n+1)^2}^{-n^2} d\rho(\lambda)$$

$$\leq \sum_{n=0}^{\infty} D_1 e^{\varepsilon(n+1)^2} e^{-\frac{n^2}{4D_2}} < \infty$$

if $\varepsilon < 1/4D_2$. $\square$

*Remark.* If in addition $q \in L^1((0,\infty))$, then the corresponding Schrödinger operator is bounded from below and hence $d\rho$ has compact support on $(-\infty, 0]$. This fact will be useful in the scattering theoretic context at the end of Section 8.

The estimate (6.8), in the case of non-Dirichlet boundary conditions at $x = 0_+$, appears to be due to Marchenko [26]. Since it is a fundamental ingredient in the inverse spectral problem, it generated considerable attention; see, for instance, [12], [18], [19], [20], [22], [27], [28, §2.4]. The case of a Dirichlet boundary at $x = 0_+$ was studied in detail by Levitan [20]. These authors, in addition to studying the spectral asymptotics of $\rho(\lambda)$ as $\lambda \downarrow -\infty$, were also particularly interested in the asymptotics of $\rho(\lambda)$ as $\lambda \uparrow \infty$ and established Theorem 6.1 (and (A.9)). In the latter context, we also refer to Bennewitz [4], Harris [16], and the literature cited therein. In contrast to these activities, we were not able to find estimates of the type (6.7) (which implies (6.8)) and (6.11) in the literature.

## 7. The relation between $A$ and $\rho$: Distributional form, II

We can now extend Theorem 5.2 to all cases.

THEOREM 7.1. *Let $f \in C_0^\infty((0,\infty))$ and suppose $b = \infty$. Assume $q$ satisfies (1.3) and let $d\rho$ be the associated spectral measure and $A$ the associated A-function. Then (5.16) and (5.17) hold.*

*Proof.* Suppose $f \in C_0^\infty((0,a))$. For $R > a$, we can find $h(R)$ so $d\rho_{R,h(R)} \underset{R\to\infty}{\to} d\rho$ (by Proposition 6.3) weakly. By Proposition 6.2 we have uniform bounds on $\int_0^\infty (1+\lambda^2)^{-1} d\rho_{R,h(R)}$ and by Theorem 6.6 on $\int_{-\infty}^0 e^{2a(\sqrt{-\lambda})} d\rho_{R,h(R)}$. Since the $\alpha$ integral in (5.15) is bounded by $C(1+\lambda^2)^{-1}$ for $\lambda > 0$ and by $Ce^{2a\sqrt{-\lambda}}$ for $\lambda \leq 0$, the right-hand side of (5.15) converges as $R \to \infty$ to the $d\rho$ integral. By Theorem 1.3, $A$ is independent of $R$ for $\alpha \in (0,a)$ and $R > a$, so the left-hand side of (5.15) is constant. Thus, (5.15) holds for $d\rho$. □

## 8. The relation between $A$ and $\rho$, III: Abelian limits

For $f \in C_0^\infty(\mathbb{R})$, define for $\lambda \in \mathbb{R}$

$$(8.1) \qquad Q(f)(\lambda) = \int_{-\infty}^\infty f(\alpha) \frac{\sin(2\alpha\sqrt{\lambda})}{\sqrt{\lambda}} d\alpha$$

and then

$$T(f) = -2\int_{\mathbb{R}} Q(f)(\lambda)\, d\rho(\lambda) \tag{8.2}$$

$$= \int_{-\infty}^{\infty} \tilde{A}(\alpha) f(\alpha)\, d\alpha. \tag{8.3}$$

We have proven in (5.16), (5.17) that for $f \in C_0^\infty(\mathbb{R})$, the two expressions (8.2), (8.3) define the same $T(f)$. We only proved this for odd $f$'s but both integrals vanish for even $f$'s. We will use (8.2) to extend to a large class of $f$, but need to exercise some care not to use (8.3) except for $f \in C_0^\infty(\mathbb{R})$.

$Q(f)$ can be defined as long as $f$ satisfies

$$|f(\alpha)| \leq C_k e^{-k|\alpha|}, \qquad \alpha \in \mathbb{R} \tag{8.4}$$

for all $k > 0$. In particular, a simple calculation shows that

$$f(\alpha) = (\pi\varepsilon)^{-\frac{1}{2}}\left[e^{-(\alpha-\alpha_0)^2/\varepsilon}\right] \Rightarrow Q(f)(\lambda) = \frac{\sin(2\alpha_0\sqrt{\lambda})}{\sqrt{\lambda}}\,e^{-\varepsilon\lambda}. \tag{8.5}$$

We use $f(\alpha,\alpha_0,\varepsilon)$ for the function $f$ in (8.5).

For $\lambda \geq 0$, repeated integrations by parts show that

$$|Q(f)(\lambda)| \leq C(1+\lambda^2)^{-1}\left[\|f\|_1 + \left\|\frac{d^3 f}{d\alpha^3}\right\|_1\right], \tag{8.6}$$

where $\|\cdot\|_1$ represents the $L^1(\mathbb{R})$-norm. Moreover, essentially by repeating the calculation that led to (8.5), we see that for $\lambda \leq 0$,

$$|Q(f)(\lambda)| \leq C e^{\varepsilon|\lambda|}\left\|e^{+\alpha^2/\varepsilon} f\right\|_\infty. \tag{8.7}$$

We conclude

PROPOSITION 8.1. *If $\int_{\mathbb{R}} d\rho(\lambda)(1+\lambda^2)^{-1} < \infty$ (always true!) and $\int_{-\infty}^0 e^{-\varepsilon_0\lambda} d\rho(\lambda) < \infty$ (see Theorem 6.7 and the remark following its proof), then using (8.2), $T(\,\cdot\,)$ can be extended to $C^3(\mathbb{R})$ $f$'s that satisfy $e^{\alpha^2/\varepsilon_0} f \in L^\infty(\mathbb{R})$ for some $\varepsilon_0 > 0$ and $\frac{d^3 f}{d\alpha^3} \in L^1(\mathbb{R})$, and moreover,*

$$|T(f)| \leq C\left[\left\|\frac{d^3 f}{d\alpha^3}\right\|_1 + \left\|e^{\alpha^2/\varepsilon_0} f\right\|_\infty\right] \tag{8.8}$$

$$:= C|||f|||_{\varepsilon_0}.$$

Next, fix $\alpha_0$ and $\varepsilon_0 > 0$ so that $\int_{-\infty}^0 e^{-\varepsilon_0\lambda} d\rho(\lambda) < \infty$. If $0 < \varepsilon < \varepsilon_0$, $f(\alpha,\alpha_0,\varepsilon)$ satisfies $|||f|||_{\varepsilon_0} < \infty$ so we can define $T(f)$. Fix $g \in C_0^\infty(\mathbb{R})$ with $g := 1$ on $(-2\alpha_0, 2\alpha_0)$. Then $|||f(\,\cdot\,,\alpha_0,\varepsilon)(1-g)|||_{\varepsilon_0} \to 0$ as $\varepsilon \downarrow 0$. So

$$\lim_{\varepsilon\downarrow 0} T(f(\,\cdot\,,\alpha_0,\varepsilon)) = \lim_{\varepsilon\downarrow 0} T(gf(\,\cdot\,,\alpha_0,\varepsilon)).$$

For $gf$, we can use the expression (8.3). $f$ is approximately $\delta(\alpha - \alpha_0)$ so standard estimates show if $\alpha_0$ is a point of Lebesgue continuity of $\tilde{A}(\alpha)$, then

$$\int_{-\infty}^{\infty} f(\alpha, \alpha_0, \varepsilon) g(\alpha) \tilde{A}(\alpha)\, d\alpha \underset{\varepsilon \downarrow 0}{\to} \tilde{A}(\alpha_0).$$

Since $A - q$ is continuous, points of Lebesgue continuity of $A$ exactly are points of Lebesgue continuity of $q$. We have therefore proven

THEOREM 8.2.  *Suppose either $b < \infty$ and $q \in L^1((0,b))$ or $b = \infty$, and then either $q \in L^1((0,\infty))$ or $q \in L^1((0,a))$ for all $a < \infty$ and*

$$q(x) \geq -Cx^2, \qquad x \geq R$$

*for some $R > 0$, $C > 0$. Let $\alpha_0 \in (0,b)$ and be a point of Lebesgue continuity of $q$. Then*

$$(8.9) \qquad A(\alpha_0) = -2 \lim_{\varepsilon \downarrow 0} \int_{\mathbb{R}} e^{-\varepsilon \lambda} \frac{\sin(2\alpha_0 \sqrt{\lambda})}{\sqrt{\lambda}}\, d\rho(\lambda).$$

We briefly illustrate the rate of convergence as $\varepsilon \downarrow 0$ in (8.9) in the special case where $q^{(0)}(x) = 0$, $x \geq 0$. Then $d\rho^{(0)}(\lambda) = \pi^{-1} \chi_{[0,\infty)}(\lambda) \sqrt{\lambda}\, d\lambda$ and formula 3.9521 of [15] (changing variables to $k = \sqrt{\lambda} \geq 0$) yield

$$(8.10) \qquad A^{(0)}(\alpha) = -2\pi^{-1} \lim_{\varepsilon \downarrow 0} \int_0^{\infty} e^{-\varepsilon \lambda} \sin(2\alpha\sqrt{\lambda})\, d\lambda$$

$$= -2\alpha \pi^{-\frac{1}{2}} \lim_{\varepsilon \downarrow 0} \varepsilon^{-\frac{3}{2}} \exp\left(-\frac{\alpha^2}{\varepsilon}\right) = 0, \qquad \alpha \geq 0.$$

Finally, we specialize (8.9) to the scattering theoretic setting. Assuming $q \in L^1((0,\infty); (1+x)\, dx)$, the corresponding Jost solution $f(x,z)$ is defined by

$$(8.11) \quad f(x,z) = e^{i\sqrt{z}\, x} - \int_x^{\infty} \frac{\sin(\sqrt{z}\,(x-x'))}{\sqrt{z}} q(x') f(x',z)\, dx', \quad \operatorname{Im}(\sqrt{z}) \geq 0$$

and the corresponding Jost function, $F(\sqrt{z})$, and scattering matrix, $S(\lambda)$, $\lambda \geq 0$, then read

$$(8.12) \qquad F(\sqrt{z}) = f(0_+, z),$$

$$(8.13) \qquad S(\lambda) = \overline{F(\sqrt{\lambda})}/F(\sqrt{\lambda}), \qquad \lambda \geq 0.$$

The spectrum of the Schrödinger operator in $L^2((0,\infty))$ associated with the differential expression $-\frac{d^2}{dx^2} + q(x)$ and a Dirichlet boundary condition at $x = 0_+$ is simple and of the type

$$\{-\kappa_j^2 < 0\}_{j \in J} \cup [0, \infty).$$

Here $J$ is a finite (possibly empty) index set, $\kappa_j > 0$, $j \in J$, and the essential spectrum is purely absolutely continuous. The corresponding spectral measure explicitly reads

$$(8.14) \qquad d\rho(\lambda) = \begin{cases} \pi^{-1}|F(\sqrt{\lambda})|^{-2}\sqrt{\lambda}\,d\lambda, & \lambda \geq 0, \\ \sum_{j\in J} c_j \delta(\lambda + \kappa_j^2)\,d\lambda, & \lambda < 0, \end{cases}$$

where

$$(8.15) \qquad c_j = \|\varphi(\,\cdot\,, -\kappa_j^2)\|_2^{-2}, \qquad j \in J$$

are the norming constants associated with the eigenvalues $\lambda_j = -\kappa_j^2 < 0$. Here $\varphi(x,z)$ (which has been introduced in (3.6a) and (3.10)) and $f(x,z)$ in (8.11) are linearly dependent precisely for $z = -\kappa_j^2$, $j \in J$.

Since

$$|F(\sqrt{\lambda})| = \prod_{j\in J}\left(1 + \frac{\kappa_j^2}{\lambda}\right)\exp\left(\frac{1}{\pi}P\int_0^\infty \frac{\delta(\lambda')\,d\lambda'}{\lambda - \lambda'}\right), \qquad \lambda \geq 0,$$

where $P\int_0^\infty$ denotes the principal value symbol and $\delta(\lambda)$ the corresponding scattering phase shift, that is, $S(\lambda) = \exp(2i\delta(\lambda))$, $\delta(\lambda) \underset{\lambda\uparrow\infty}{\to} 0$, the scattering data

$$\{-\kappa_j^2, c_j\}_{j\in J} \cup \{S(\lambda)\}_{\lambda \geq 0}$$

uniquely determine the spectral measure (8.14) and hence $A(\alpha)$. Inserting (8.14) into (8.9) then yields the following expression for $A(\alpha)$ in terms of scattering data.

THEOREM 8.3. *Suppose that $q \in L^1((0, \infty); (1 + x)\,dx)$. Then*

$$(8.16) \qquad A(\alpha) = -2\sum_{j\in J} c_j \kappa_j^{-1}\sinh(2\alpha\kappa_j)$$

$$- 2\pi^{-1}\lim_{\varepsilon\downarrow 0}\int_0^\infty e^{-\varepsilon\lambda}|F(\sqrt{\lambda})|^{-2}\sin(2\alpha\sqrt{\lambda})\,d\lambda$$

*at points $\alpha \geq 0$ of Lebesgue continuity of $q$.*

*Remark.* In great generality $|F(k)| \to 1$ as $k \to \infty$, so one cannot take the limit in $\varepsilon$ inside the integral in (8.16). In general, though, one can can replace $|F(\sqrt{\lambda})|^{-2}$ by $(|F(\sqrt{\lambda})|^{-2} - 1) \equiv X(\lambda)$ and ask if one can take a limit there. As long as $q$ is $C^2((0,\infty))$ with $q'' \in L^1((0,\infty))$, it is not hard to see that as $\lambda \to \infty$

$$X(\lambda) = -\frac{q(0)}{2\lambda} + O(\lambda^{-2}).$$

Thus, if $q(0) = 0$, then

(8.17) $$A(\alpha) = -2 \sum_{j \in J} c_j \kappa_j^{-1} \sinh(2\alpha \kappa_j)$$
$$- 2\pi^{-1} \int_0^\infty (|F(\sqrt{\lambda})|^{-2} - 1) \sin(2\alpha \sqrt{\lambda}) \, d\lambda.$$

The integral in (8.17) is only conditionally convergent if $q(0) \neq 0$.

We note that in the present case where $q \in L^1((0, \infty); (1 + x) \, dx)$, the representation (1.17) of the $m$-function in terms of the $A(\alpha)$-amplitude was considered in a paper by Ramm [31] (see also [32, pp. 288–291]).

## 9. The relation between $A$ and $\rho$, IV: Remarks

Here is a totally formal way of understanding why (5.1) is true. We start with the basic representation without errors,

(9.1) $$m(-\kappa^2) = m_0(-\kappa^2) - \int_0^\infty A(\alpha) e^{-2\alpha\kappa} \, d\alpha.$$

Pretend we can analytically continue from $\kappa$ real to $\kappa = -ik$ (at which point $-\kappa^2$ is $k^2 + i0$). Then

(9.2) $$m(k^2 + i0) = m_0(k^2 + i0) - \int_0^\infty A(\alpha) e^{2i\alpha k} \, d\alpha.$$

This normally cannot be literally true. In many cases, $A(\alpha) \to \infty$ at infinity (although for the case $q(x) = \text{constant} > 0$, which we discuss later, it is true). But this is only a formal argument.

Taking imaginary parts and using for $\alpha, \alpha_0 > 0$ that

(9.3) $$\int_0^\infty \sin(2\alpha k) \sin(2\alpha_0 k) \, dk = \frac{2\pi}{8} \delta(\alpha - \alpha_0)$$

(which follows from $\int_{-\infty}^\infty e^{i\alpha k} \, dk = 2\pi\delta(\alpha)$), we conclude that for $\alpha_0 > 0$,

(9.4) $$A(\alpha_0) = -\frac{4}{\pi} \int_0^\infty \sin(2\alpha_0 k) \text{Im}(m_0(k^2 + i0)) \, dk$$
$$= -2 \int_0^\infty \sin(2\alpha_0 \sqrt{\lambda}) \left[ \frac{\text{Im}(m_0(\lambda + i0))}{\pi} \right] \frac{d\lambda}{\sqrt{\lambda}},$$

which, given (1.12), is just (5.1).

As explained in [35], a motivation for $A$ is the analogy to the $m$-function for a tridiagonal Jacobi matrix. For this point of view, the relation (5.1) is an important missing link. The analog of (1.7) in the discrete case is

(9.5) $$m(z) = -\sum_{n=0}^\infty \frac{\gamma_n}{z^{n+1}}.$$

The coefficients of $\gamma_n$ of the Taylor series at infinity are the analog of $A(\alpha)$. In this case, the spectral measure is finite and of finite support (if the Jacobi matrix is bounded) and

$$(9.6) \qquad m(z) = \int_{\mathbb{R}} \frac{d\rho(\lambda)}{\lambda - z}$$

so that (9.5) implies that

$$(9.7) \qquad \gamma_n = \int_{\mathbb{R}} \lambda^n \, d\rho(\lambda).$$

(5.1) should be then thought of as the analog of (9.7) for the continuum case.

Perhaps the most important consequence of (9.7) is the implied positivity condition of the $\gamma$'s — explicitly, that

$$\sum_{m,n=0}^{N} \gamma_{n+m} a_m \overline{a_n} \geq 0$$

for all $(a_0, \ldots, a_N) \in \mathbb{C}^{N+1}$.

Recall (see, e.g., Gel'fand-Vilenkin [13, §II.5]) that Krein proved the following fact:

THEOREM 9.1.  *A continuous even function $f$ on $\mathbb{R}$ has the property that*

$$(9.8) \qquad \int_{\mathbb{R}^2} f(x-y) \varphi(y) \overline{\varphi(x)} \, dx dy \geq 0$$

*for all even functions $\varphi \in C_0^\infty(\mathbb{R})$ if and only if there are finite positive measures $d\mu_1$ and $d\mu_2$ on $[0, \infty)$ so that $\int_0^\infty e^{a\lambda} \, d\mu_2(\lambda) < \infty$ for all $a > 0$ and so that*

$$(9.9) \qquad f(x) = \int_0^\infty \cos(\lambda x) \, d\mu_1(\lambda) + \int_0^\infty \cosh(\lambda x) \, d\mu_2(\lambda).$$

Using the extension in Gel'fand-Vilenkin to distributional $f$ (cf. [13, §II.6.3, Theorem 5]), one obtains

THEOREM 9.2.  *Let $\tilde{A}(\alpha)$ be the distribution of Theorem 5.3. Let $B(\alpha) = -\tilde{A}'(\alpha)$ be the distributional derivative of $\tilde{A}$. Then*

$$(9.10) \qquad \int_{\mathbb{R}} B(\alpha) \left[ \int_{\mathbb{R}} \varphi(\beta) \overline{\varphi(\beta - \alpha)} \, d\beta \right] d\alpha \geq 0$$

*for all even $\varphi \in C_0^\infty(\mathbb{R})$. Moreover, if $\tilde{A}$ is a distribution related to a signed measure, $d\rho$, by (5.16), then (9.10) is equivalent to the positivity of the measure $d\rho$.*

As discussed in [13], the measures $d\mu_j$ in (9.9) may not be unique. Our theory illuminates this fact. If $q$ is in the limit circle case at infinity, then distinct boundary conditions lead to distinct spectral measures but the same $A$-function, so the same $\tilde{A}$ and the same $B = \tilde{A}'$. Thus, we have additional examples of nonuniqueness. The growth restrictions on $f$ which guarantee uniqueness in (9.9) (e.g., $\int_\mathbb{R} e^{-cx^2} f(x)\,dx < \infty$ for all $c > 0$) are not unrelated to the standard $q(x) \geq -Cx^2$ that leads to the limit point case at $\infty$ for the Schrödinger differential expression $-\frac{d^2}{dx^2} + q(x)$.

Next we turn to the relation between $A$ and the Gel'fand-Levitan transformation kernel $L$ in [12]. For the function $L(x,y)$ associated to Dirichlet boundary conditions at $x = 0$, satisfying (cf. (3.6a), (3.10))

$$\frac{\sin(\sqrt{z}\,x)}{\sqrt{z}} = \varphi(x,z) + \int_0^x L(x,x')\varphi(x',z)\,dx',$$

we claim that

$$(9.11) \qquad A(\alpha) = -2 \frac{\partial}{\partial y} L(2\alpha, y)\bigg|_{y=0_+}.$$

We will first proceed formally without worrying about regularity conditions. Detailed discussions of transformation operators can be found, for instance, in [11], [21, Ch. 1], [22], [24], [25], [26], [28, Ch. 1], [30, Ch. VIII], [36], [37], and, in the particular case of scattering theory, in [2, Chs. I and V], [8], and [29]. Let $d\rho(\lambda)$ be the spectral measure for $-\frac{d^2}{dx^2} + q(x)$ and

$$(9.12) \qquad d\rho^{(0)}(\lambda) = \pi^{-1}\chi_{[0,\infty)}(\lambda)\sqrt{\lambda}\,d\lambda$$

the spectral measure for $-\frac{d^2}{dx^2}$ (both corresponding to the Dirichlet boundary condition parameter $h = \infty$ at $x = 0$), and define $d\sigma = d\rho - d\rho^{(0)}$. Then $L$ is defined as follows [12]. Let

$$(9.13) \qquad F(x,y) = \int_{-\infty}^{\infty} \frac{[1-\cos(\sqrt{\lambda}\,x)][1-\cos(\sqrt{\lambda}\,y)]}{\lambda^2}\,d\sigma(\lambda)$$

and

$$(9.14) \qquad k(x,y) = \frac{\partial^2 F}{\partial x \partial y}(x,y) \text{``$=$''} \int_{-\infty}^{\infty} \frac{\sin(\sqrt{\lambda}\,x)\sin(\sqrt{\lambda}\,y)}{\lambda}\,d\sigma(\lambda),$$

where the final "$=$" is formal since the integral may not converge absolutely. $L$ satisfies the following nonlinear Gel'fand-Levitan equation,

$$(9.15) \qquad k(x,y) = L(x,y) + \int_0^y L(x,x')L(y,x')\,dx',$$

$$(9.16) \qquad L(x, 0_+) = 0, \qquad L(x,x) = -\frac{1}{2}\int_0^x q(x')\,dx'.$$

Thus, formally by (9.15) and (9.16),

$$(9.17) \qquad \frac{\partial k}{\partial y}(x,y)\bigg|_{y=0_+} = \frac{\partial L}{\partial y}(x,y)\bigg|_{y=0_+},$$

and then by (9.14)

$$\left.\frac{\partial k}{\partial y}(x,y)\right|_{y=0_+} \text{``=''} \int_{-\infty}^{\infty} \frac{\sin(\sqrt{\lambda}\,x)}{\sqrt{\lambda}}\,d\sigma(\lambda), \tag{9.18}$$

which, by (5.1), says that (9.11) holds.

Alternatively, one can derive (9.11) as follows. Suppose $Q \in L^1((0,\infty))$ coincides with $q$ on the interval $[0,\alpha]$, is real-valued, and of compact support. Denote by $f_Q(x,z)$, $F_Q(\sqrt{z})$, and $L_Q(x,x')$ the Jost solution, Jost function, and transformation kernel (satisfying (9.15), (9.16)) associated with $Q$. Then (cf. [5, §V.2]),

$$\frac{f_Q(x,z)}{F_Q(\sqrt{z})} = e^{i\sqrt{z}\,x} + \int_x^{\infty} L_Q(x',x) e^{i\sqrt{z}\,x'}\,dx', \tag{9.19}$$

and

$$u_Q(x,z) = \frac{f_Q(x,z)}{F_Q(\sqrt{z})}, \qquad u_Q(0_+,z) = 1, \tag{9.20a}$$

$$u_Q(\cdot,z) \in L^2((0,\infty)), \qquad z \in \mathbb{C}\backslash\mathbb{R} \tag{9.20b}$$

is the unique Weyl solution association with $Q$. Thus, the normalization of $u_Q$ in (9.20a), (9.19), $L_Q(0_+,0_+) = 0$, and (1.7) then yield

$$m_Q(z) = u_Q'(0_+,z) = i\sqrt{z} + \int_0^{\infty} \left(\left.\frac{\partial}{\partial x}L_Q(x',x)\right|_{x=0_+}\right) e^{i\sqrt{z}\,x'}\,dx'. \tag{9.21}$$

Identifying $z = -\kappa^2$, $x' = 2\alpha$, a comparison with (1.17) then implies

$$A_Q(\alpha) = -2L_{Q,y}(2\alpha,0_+).$$

Since by Theorem 1.3 and the following remark, $A(\alpha)$ only depends on $q(x) = Q(x)$ for $x \in [0,\alpha]$, and $L(x,y)$ depends on $q(x') = Q(x')$ for $x' \in [0,(x+y)/2]$ with $0 \le y \le x \le 2\alpha$ (cf. [5, eq. (III.1.11)], [28, pp. 19, 20]), one concludes (9.11).

Next, we want to note that (5.1) sometimes does not represent a conditionally convergent integral; that is,

$$A(\alpha) = -2 \lim_{R\to\infty} \int_{-\infty}^{R} \lambda^{-\frac{1}{2}} \sin(2\alpha\sqrt{\lambda})\,d\rho(\lambda) \tag{9.22}$$

can fail. Indeed, it even fails in the case $b < \infty$, $h = \infty$, and $q^{(0)}(x) = 0$, $0 \le x \le b$. For in that case (see (4.4)),

$$m^{(0)}(-\kappa^2) = -\frac{\kappa + \kappa e^{-2\kappa b}}{1 - e^{-2\kappa b}}.$$

Straightforward residue calculus then implies that

$$d\rho^{(0)}(\lambda) = \sum_{n=1}^{\infty} w_n \delta(E - E_n), \tag{9.23a}$$

with

$$E_n = \frac{\pi^2 n^2}{b^2} \tag{9.23b}$$

and

$$w_n = \frac{2\pi^2 n^2}{b^3} \tag{9.23c}$$

(the reader might want to check that this is consistent with $\int_0^R d\rho(\lambda) \underset{R\to\infty}{\sim} \frac{2}{3\pi} R^{3/2}$).

Thus,

$$\int_{-\infty}^{R} \lambda^{-\frac{1}{2}} \sin(2\alpha\sqrt{\lambda}) \, d\rho^{(0)}(\lambda) = \sum_{n \leq bR^{1/2}/\pi} \frac{2\pi n}{b^2} \sin(2\pi\alpha n/b)$$

is not conditionally convergent as $R \to \infty$.

Given the known asymptotics for the eigenvalues and weights when $b < \infty$ (cf., e.g., [23, §1.2]), one can see that (9.22) never holds if $b < \infty$. There are also cases with $b = \infty$, where it is easy to see the integral cannot be conditionally convergent. If

$$q(x) = x^\beta, \qquad \beta > 0$$

then WKB analysis (see, e.g., [38, §7.1]) shows that

$$E_n \underset{n\uparrow\infty}{=} [Cn + O(1)]^\gamma,$$

where $\gamma^{-1} = \frac{1}{2} + \beta^{-1}$ and $w_n = CE_n^{1-1/\beta}(1 + o(1))$. As long as $\beta > 2$, $w_n E_n^{-\frac{1}{2}} \underset{n\to\infty}{\to} \infty$, and so the integral is not conditionally convergent.

Another canonical scenario displaying this phenomenon is provided by the scattering theoretic setting discussed at the end of Section 8. In fact, assuming $q \in L^1((0,\infty); (1+x)\,dx)$, one sees that

$$|F(k)| \underset{k\uparrow\infty}{=} 1 + o(k^{-1}) \tag{9.24}$$

(cf. [5, eq. II.4.13] and apply the Riemann-Lebesgue lemma. Actually, one only needs $q \in L^1((0,\infty))$ for the asymptotic results on $F(k)$ as $k \uparrow \infty$ but we will ignore this refinement in the following.) A comparison of (9.24) and (8.16) then clearly demonstrates the necessity of an Abelian limit in (8.16). Even replacing $d\rho$ in (8.9) by $d\sigma = d\rho - d\rho^{(0)}$ (cf. (8.10)); that is, effectively replacing $|F(\sqrt{\lambda})|^{-2}$ by $[|F(\sqrt{\lambda})|^{-2} - 1]$ in (8.16) still does not necessarily produce an absolutely convergent integral in (8.16).

The latter situation changes upon increasing the smoothness properties of $q$ since, for example, assuming $q \in L^1((0,\infty);(1+x)\,dx)$, $q' \in L^1((0,\infty))$, yields

$$|F(k)|^{-2} - 1 \underset{k\uparrow\infty}{=} O(k^{-2})$$

as detailed high-energy considerations (cf. [14]) reveal. Indeed as we saw at the end of Section 8, if $q'' \in L^1((0,\infty))$, then the integral one gets is absolutely convergent if and only if $q(0) = 0$.

Unlike the oscillator-like cases, though, the integrals in the scattering theory case are conditionally convergent.

These examples allow us to say something about the following question raised by R. del Rio [9]. Does $|m(-\kappa^2) + \kappa|$ stay bounded as $z = -\kappa^2$ moves along the curve $\mathrm{Im}(z) = a_0 > 0$ with $\mathrm{Re}(z) \to \infty$? In general, the answer is no. The $m(-\kappa^2)$ of (4.4) has

$$G(-\kappa^2) := |m(-\kappa^2) + \kappa| = 2|\kappa|\,|e^{-2\kappa b}|/|1 - e^{-2\kappa b}|$$

so

$$G((\tfrac{\pi n}{b})^2 + ia_0)/\sqrt{|E_n|} \underset{n\to\infty}{\to} \infty$$

($E_n = (\pi n/b)^2$, $n \in \mathbb{N}$ denoting the corresponding eigenvalues) showing $|m(-\kappa^2)|$ is not even bounded by $C|\kappa|$ on the curve. Similarly, in the case $q(x) = x^\beta$, one infers that $|m(-\kappa^2) + \kappa|/|\kappa|$ is unbounded at $E_n + ia_0$ as long as $\beta > 2$.

As a final issue related to the representation (5.1), we discuss the issue of bounds on $A$ when $|q(x)| \leq Cx^2$. We have two general bounds on $A$: the estimate of [35] (see ((1.16)),

$$(9.25) \quad |A(\alpha) - q(\alpha)| \leq \left[\int_0^\alpha |q(y)|\,dy\right]^2 \exp\left[\alpha \int_0^\alpha |q(y)|\,dy\right],$$

and the estimate we will prove in the next section (Theorem 10.2),

$$(9.26) \quad |A(\alpha)| \leq \frac{\gamma(\alpha)}{\alpha} I_1(2\alpha\gamma(\alpha)),$$

where $|\gamma(\alpha)| = \sup_{0\leq x\leq \alpha} |q(x)|^{1/2}$ and $I_1(\cdot)$ is the modified Bessel function of order one (cf., e.g., [1, Ch. 9]). Since ([1, p. 375])

$$(9.27) \quad 0 \leq I_1(x) \leq e^x, \qquad x \geq 0,$$

we conclude that

$$(9.28) \quad |A(\alpha)| \leq \sqrt{C}\,e^{2\sqrt{C}\,\alpha^2}$$

if $|q(x)| \leq Cx^2$. This is a pointwise bound related to the integral bounds on $A(\alpha)$ implicit in Lemma 6.5.

## 10. Examples, I: Constant $q$

We begin with the case $b = \infty$, $q(x) = q_0$, $x \geq 0$, with $q_0$ a real constant. We claim

THEOREM 10.1. *If $b = \infty$ and $q(x) = q_0$, $x \geq 0$, then if $q_0 > 0$,*

$$A(\alpha) = \frac{q_0^{\frac{1}{2}}}{\alpha} J_1(2\alpha q_0^{\frac{1}{2}}), \tag{10.1}$$

*where $J_1(\,\cdot\,)$ is the Bessel function of order one (cf., e.g., [1, Ch. 9]); and if $q_0 < 0$,*

$$A(\alpha) = \frac{(-q_0)^{\frac{1}{2}}}{\alpha} I_1(2\alpha(-q_0)^{\frac{1}{2}}), \tag{10.2}$$

*with $I_1(\,\cdot\,)$ the corresponding modified Bessel function.*

*Proof.* We use the following formula ([15, 6.6233]),

$$\int_0^\infty e^{-ax} J_1(bx) \frac{dx}{x} = \frac{\sqrt{a^2+b^2} - a}{b}, \quad a > 0,\ b \in \mathbb{R}, \tag{10.3}$$

$$\int_0^\infty e^{-ax} I_1(bx) \frac{dx}{x} = \frac{\sqrt{a^2-b^2} - a}{b}, \quad a > 0,\ |b| < a,\ b \in \mathbb{R}. \tag{10.4}$$

From this we see that

$$-\kappa - \int_0^\infty e^{-2\alpha\kappa} \frac{q_0^{\frac{1}{2}}}{\alpha} J_1(2\alpha q_0^{\frac{1}{2}})\, d\alpha = -\frac{\sqrt{(2\kappa)^2 + (2q_0^{\frac{1}{2}})^2} - 2\kappa}{2} - \kappa$$

$$= -\sqrt{\kappa^2 + q_0},$$

which is the $m$-function for $b = \infty$, $q(x) = q_0$, $x \geq 0$. By the uniqueness of inverse Laplace transforms, this proves (10.1) and incidentally a formula like (1.17) without error term holds. The argument for (10.2) using (10.4) is similar. □

*Remarks.* 1. This suggests that a formula like (1.17) holds if $q$ is bounded. We will prove that below (see Theorem 10.3).

2. Our original derivation of the formula used (5.1), the known formula for $d\rho(\lambda)$, and an orgy of Besselology.

This example is especially important because of a monotonicity property:

THEOREM 10.2. *Let $|q_1(x)| \leq -q_2(x)$ on $[0,a]$ with $a \leq \min(b_1, b_2)$. Then $|A_1(\alpha)| \leq -A_2(\alpha)$ on $[0,a]$. In particular, for any $q$ satisfying $\sup_{0 \leq x \leq \alpha} |q(x)| < \infty$, we have that*

$$|A(\alpha)| \leq \frac{\gamma(\alpha)}{\alpha} I_1(2\alpha\gamma(\alpha)), \tag{10.5}$$

*where*

(10.6) $$\gamma(\alpha) = \sup_{0 \le x \le \alpha} (|q(x)|^{\frac{1}{2}}).$$

*In particular*, (9.27) *implies*

(10.7) $$|A(\alpha)| \le \alpha^{-1}\gamma(\alpha)e^{2\alpha\gamma(\alpha)},$$

*and if q is bounded,*

(10.8) $$|A(\alpha)| \le \alpha^{-1}\|q\|_\infty^{\frac{1}{2}} \exp(2\alpha\|q\|_\infty^{\frac{1}{2}}).$$

*Proof.* Since $A(\alpha)$ is only a function of $q$ on $[0, \alpha)$, we can suppose that $b_1 = b_2 = \infty$ and $q_1(x) = q_2(x) = 0$ for $x > a$. By a limiting argument, we can suppose that $q_j$ are $C^\infty([0, a])$. We can then use the expansion of [35, §2],

(10.9) $$-A(\alpha) = -q(\alpha) + \sum_{n=2}^{\infty} \frac{(-1)^n}{2^{n-2}} \int_{R_n(\alpha)} q(x_1) \ldots q(x_n)$$
$$\times\ dx_1 \ldots dx_n\, d\ell_1 \ldots d\ell_{n-2},$$

where $R_n(\alpha)$ is a complicated region on $\{x, \ell\}$ space that is $q$ independent (given by (2.19) from [35]). The monotonicity result follows immediately from this expression. (10.5) then follows from (10.2), and (10.7), (10.8) from (9.27). □

*Remarks.* 1. The expansion

$$I_1(x) = \sum_{n=0}^{\infty} \frac{(\frac{1}{2}z)^{2n+1}}{n!(n+1)!}$$

allows one to compute exactly the volume of the region $R_n(\alpha)$ of [35], *viz.*,

$$|R_n(\alpha)| = \frac{\alpha^{2n-2}}{(n-1)!n!}.$$

The bounds in [35] only imply $|R_n(\alpha)| \le \frac{\alpha^{2n-2}}{(n-2)!}$ and are much worse than the actual answer for large $n$!

2. For $\alpha$ small, (10.7) is a poor estimate and one should use (9.25) which implies that $|A(\alpha) \le \|q\|_\infty + \alpha^2\|q\|_\infty^2 e^{\alpha^2\|q\|_\infty}$.

This lets us prove

THEOREM 10.3.   *Let $h = \infty$ and $q \in L^\infty((0, \infty))$. Suppose $\kappa^2 > \|q\|_\infty$. Then*

(10.10) $$m(-\kappa^2) = -\kappa - \int_0^\infty A(\alpha)e^{-2\alpha\kappa}\, d\alpha$$

(*with a convergent integral and no error term*).

*Proof.* Let $q_n = q\chi_{[0,n]}(x)$. Let $m_n, A_n$ be the $m$-function and $A$-amplitude, respectively, for $q_n$. Then

(1) $m_n(z) \to m(z)$ for $z \in \mathbb{C} \setminus [-\|q\|_\infty, \infty)$.

(2) $A_n(\alpha) \to A(\alpha)$ pointwise (since $A_n(\alpha) = A(\alpha)$ if $n > \alpha$).

(3) (10.10) holds for $q_n$ since $q_n \in L^1((0, \infty))$ (see Theorem 1.2).

(4) $|A_n(\alpha)| \leq \alpha^{-1} \|q\|_\infty^{\frac{1}{2}} \exp(2\alpha \|q\|_\infty^{\frac{1}{2}})$. This is (10.8).

(5) $|A_n(\alpha)| \leq \|q\|_\infty [1 + \alpha^2 \|q\|_\infty \exp(\alpha^2 \|q\|_\infty)]$. This is (9.25).

The dominated convergence theorem thus implies that (10.10) holds for $q \in L^\infty((0, \infty))$. □

*Remarks.* 1. If $\inf \operatorname{supp}(d\rho) = -E_0$ with $E_0 > 0$, then $m(z)$ has a singularity at $z = -E_0$ so we cannot expect that $|A(\alpha)| \leq Ce^{2(E_0-\varepsilon)\alpha}$ for any $\varepsilon > 0$. Thus, $A$ must grow exponentially as $\alpha \to \infty$. One might naively guess that if $\inf \operatorname{supp}(d\rho) = E_0$ with $E_0 > 0$, then $A(\alpha)$ decays exponentially, but this is false in general. For example, if $q(x) = q_0 > 0$, then by (10.1) for $\alpha$ large, $A(\alpha) \sim -\pi^{-\frac{1}{2}} q_0^{\frac{1}{4}} \alpha^{-\frac{3}{2}} \cos(2q_0^{\frac{1}{2}} \alpha + \frac{\pi}{4}) + O(\alpha^{-2})$ by the known asymptotics of $J_1$ ([1, p. 364]).

2. For $q(x) = q_0 > 0$, $A(\alpha) \to 0$ as $\alpha \to \infty$. This leads one to ask if perhaps $A(\alpha) \to 0$ for all cases where $\operatorname{supp}(\rho) \subset [0, \infty)$ or at least if $q(x) > 0$. It would be interesting to know the answer even for the harmonic oscillator.

3. We have proven exponential bounds on $A(\alpha)$ as $\alpha \to \infty$ for the cases $q \in L^1((0, \infty))$ and $q \in L^\infty((0, \infty))$, but not even for $L^1((0, \infty)) + L^\infty((0, \infty))$. One might guess that $\sup_{x>0}(\int_x^{x+1} |q(y)| \, dy)$ suffices for such a bound.

## 11. Examples, II: Bargmann potentials

Our second set of examples involves Bargmann potentials (cf., e.g., [5, §§IV.3 and VI.1]), that is, potentials $q \in L^1((0, \infty); (1+x) \, dx)$ such that the associated Jost function $F(k)$ (cf. (8.12)) is a rational function of $k$. We explicitly discuss two simple examples and then hint at the general case.

*Case 1.* $F(k) = (k - i\kappa_1)/(k + i\kappa_1)$. Thus, $d\rho(\lambda) = d\rho^{(0)}(\lambda)$ on $[0, \infty)$ and there is a single eigenvalue at energy $\lambda = -\kappa_1^2$. There is a single norming constant, $c_1$, and it is known (cf. [5, §VI.1]) that

$$(11.1) \quad q(x) = -2 \frac{d^2}{dx^2} \ln \left[ 1 + \frac{c_1}{\kappa_1^2} \int_0^x \sinh^2(\kappa_1 y) \, dy \right].$$

In (5.1), the $\lambda > 0$ contribution to $A(\alpha)$ is the same as in the free case, and so it yields zero contribution to $A$ (cf. (8.10)). Thus,

$$A(\alpha) = -2c_1 \int_{-\infty}^0 |\lambda|^{-\frac{1}{2}} \sinh(2\alpha \sqrt{|\lambda|}) \delta(\lambda + \kappa_1^2) \, d\lambda$$

and hence

$$A(\alpha) = -\frac{2c_1}{\kappa_1} \sinh(2\alpha\kappa_1). \tag{11.2}$$

Note that $q(0+) = A(0+) = 0$ (verifying $q(0+) = A(0+)$).

*Case* 2. $F(k) = (k+i\gamma)/(k+i\beta)$, $\beta > 0$, $\gamma \geq 0$. It is known (cf. [5, p. 87]) that

$$q(x) = -8\beta^2 \left(\frac{\beta-\gamma}{\beta+\gamma}\right) \frac{e^{-2\beta x}}{(1+(\frac{\beta-\gamma}{\beta+\gamma})e^{-2\beta x})^2}. \tag{11.3}$$

The case $\gamma = 0$ corresponds to $q(x) = -2\beta^2/\cosh^2(\beta x)$ (the one-soliton potential on its odd subspace).

We claim that

$$m(-\kappa^2) = -\kappa - \frac{\gamma^2 - \beta^2}{\kappa + \gamma}, \tag{11.4}$$

for clearly, $m(-\kappa^2)$ is analytic in $\mathbb{C}\setminus[0,\infty)$ and satisfies $m(-\kappa^2) = -\kappa + O(\kappa^{-1})$ and at $\kappa = -ik$ (i.e., $E = -\kappa^2 = k^2 + i0$),

$$\mathrm{Im}(m(k^2+i0)) = \left[k - \frac{\gamma^2-\beta^2}{\gamma^2+k^2}k\right] = k\left[\frac{k^2+\beta^2}{k^2+\gamma^2}\right] = \frac{k}{|F(k)|^2}$$

consistent with (8.14). Thus, uniqueness of $m$ given $d\rho$ and the asymptotics proves (11.4). Since

$$\frac{1}{\kappa + \gamma} = 2\int_0^\infty e^{-2\alpha\kappa}e^{-2\alpha\gamma}\, d\alpha,$$

(11.4) and uniqueness of the inverse Laplace transform implies that

$$A(\alpha) = 2(\gamma^2 - \beta^2)e^{-2\alpha\gamma}. \tag{11.5}$$

Notice that $q(0+) = A(0+) = 2(\gamma^2 - \beta^2)$ and the odd soliton ($\gamma = 0$) corresponds to $A(\alpha) = -2\beta^2$, a constant.

*Remark.* Thus, we see that $A(\alpha)$ equal to a negative constant is a valid $A$-function. However, $A(\alpha)$ a positive constant, say, $A_0 > 0$, is not since then $\mathrm{Im}(m(k+i0)) = k - A_0/k$ is negative for $k > 0$ small.

In the case of a general Bargmann-type potential $q(x)$, one considers a Jost function of the form

$$F(k) = \left[\prod_{j\in J_e}\left(\frac{k-i\kappa_j}{k+i\kappa_j}\right)\right]\frac{k}{k+i\beta_0}\left[\prod_{\ell\in J_r}\left(\frac{k+i\gamma_\ell}{k+i\beta_\ell}\right)\right], \tag{11.6}$$

$$\kappa_j > 0,\ j \in J_e,\ \beta_0 \geq 0,\ \beta_\ell > 0,\ \gamma_\ell > 0,\ \ell \in J_r,$$

$$\beta_\ell \neq \gamma_{\ell'},\ \text{for all}\ \ell, \ell' \in J_r,\ \gamma_\ell \neq \kappa_j\ \text{for all}\ \ell \in J_r,\ j \in J_e,$$

with $J_e$ (resp. $J_r$) a finite (possibly empty) index set associated with the eigenvalues $\lambda_j = -\kappa_j^2 < 0$ (resp. resonances) of $q$, and $\beta_0 \geq 0$ associated with a possible zero-energy resonance of $q$. Attaching norming constants $c_j > 0$ to the eigenvalues $\lambda_j = -\kappa_j^2$, $j \in J_e$ of $q$, one then obtains

$$(11.7) \quad d\rho(\lambda) - d\rho^{(0)}(\lambda) = \begin{cases} \pi^{-1}[|F(\sqrt{\lambda})|^{-2} - 1]\sqrt{\lambda}\,d\lambda, & \lambda \geq 0 \\ \sum_{j \in J_e} c_j \delta(\lambda + \kappa_j^2)\,d\lambda, & \lambda < 0 \end{cases}$$

$$= \begin{cases} \pi^{-1} \sum_{\ell \in J_r \cup \{0\}} A_\ell(\lambda + \gamma_\ell^2)^{-1}\sqrt{\lambda}\,d\lambda, & \lambda \geq 0 \\ \sum_{j \in J_e} c_j \delta(\lambda + \kappa_j^2)\,d\lambda, & \lambda < 0. \end{cases}$$

Here $d\rho^{(0)}$ denotes the spectral measure (9.12) for the free case $q^{(0)}(x) = 0$, $x \geq 0$, and $\gamma_0 = 0$,

$$A_\ell = \begin{cases} \prod_{m \in J_r \cup \{0\}}(\beta_m^2 - \gamma_\ell^2) \prod_{\substack{n \in J_r \cup \{0\} \\ n \neq \ell}} (\gamma_n^2 - \gamma_\ell^2)^{-1}, & \ell \in J_r \cup \{0\},\ \beta_0 > 0, \\ \prod_{m \in J_r}(\beta_m^2 - \gamma_\ell^2) \prod_{\substack{n \in J_r \\ n \neq \ell}} (\gamma_n^2 - \gamma_\ell^2)^{-1}, & \ell \in J_r,\ \beta_0 = 0, \end{cases}$$

and $A_0 = 0$ if $\beta_0 = 0$.

Next, observing the spectral representation for the free Green's function associated with $q^{(0)}(x) = 0$, $x \geq 0$ and a Dirichlet boundary condition at $x = 0_+$, one computes

$$(11.8) \quad \frac{1}{\pi} \int_0^\infty \frac{\sin(\sqrt{\lambda}\,x)}{\sqrt{\lambda}} \frac{\sin(\sqrt{\lambda}\,y)}{\sqrt{\lambda}} \frac{1}{\lambda + \gamma^2} \sqrt{\lambda}\,d\lambda = e^{-\gamma x} \frac{\sinh(\gamma y)}{\gamma}, \quad x \geq y.$$

Taking into account $A^{(0)}(\alpha) = 0$, $\alpha \geq 0$ according to (8.10) (hence subtracting $d\rho^{(0)}$ in (11.7) will have no effect on computing $A(\alpha)$ using (8.16)), the $y$-derivative of the integral (11.8) at $y = 0_+$ combined with an Abelian limit $\varepsilon \downarrow 0$ yields precisely the prototype of integral (*viz.*, $\lim_{\varepsilon \downarrow 0_+} \int_0^\infty e^{-\varepsilon\lambda} \sin(2\alpha\sqrt{\lambda})(\lambda + \gamma^2)^{-1}\,d\lambda$) needed to compute $A(\alpha)$ upon inserting (11.7) into (8.16). The net result then becomes

$$(11.9) \quad A(\alpha) = -2 \sum_{j \in J_e} c_j \kappa_j^{-1} \sinh(2\alpha\kappa_j) - 2 \sum_{\ell \in J_r \cup \{0\}} A_\ell e^{-2\alpha\gamma_\ell}.$$

The corresponding potential $q(x)$ can be computed along the lines indicated in [5, Ch. IV] and is known to be continuous on $[0, \infty)$. Hence (11.9) holds for all $\alpha \geq 0$. More precisely, the condition $q \in L^1((0, \infty); (1 + x)\,dx)$ imposes certain restrictions on the possible choice of $\beta_\ell > 0$, $\gamma_\ell > 0$ in (11.6) in order to avoid isolated singularities of the type $2(x - x_0)^{-2}$ in $q(x)$. Away from such isolated singularities, (11.7) inserted into the Gel'fand-Levitan equation yields a $C^\infty$ potential $q$ (in fact, a rational function of certain exponential functions and their $x$-derivatives) upon solving the resulting linear algebraic system of equations. In particular, one obtains $q(0_+) = A(0_+) = -2 \sum_{\ell \in J_r \cup \{0\}} A_\ell$.

### Appendix A. The $B_h$ function

Throughout this paper, we have discussed the principal $m$-function, $m(z)$ given by (1.7). This is naturally associated to Dirichlet boundary conditions because the spectral measure $d\rho$ of (1.10) is a spectral measure for an operator $H$, with $u(0_+) = 0$ boundary conditions. For $h \in \mathbb{R}$, there are subsidiary $m$-functions, $m_h(z)$, associated to

$$(A.1) \qquad u'(0_+) + hu(0_+) = 0$$

boundary conditions. Our goal in this section is to present Laplace transform asymptotics for $m_h(z)$.

One defines $m_h(z)$ by

$$(A.2) \qquad m_h(z) = [hm(z) - 1]/[m(z) + h].$$

That this is associated to the boundary (A.1) is hinted at by the fact that $m(z) + h = 0$ if and only if $u'(0_+, z) + hu(0_+, z) = 0$. The function

$$F_h(\zeta) = \frac{h\zeta - 1}{\zeta + h}$$

satisfies

(i) $F_h : \mathbb{C}_+ \to \mathbb{C}_+$, where $\mathbb{C}_+ = \{z \in \mathbb{C} \mid \operatorname{Im}(z) > 0\}$,

(ii) $F_h(\zeta) = h - \frac{(1+h^2)}{\zeta+h}$,

(iii) $F_h(\zeta) - F_h(\zeta_0) = (\zeta - \zeta_0)(\zeta + h)^{-1}(\zeta_0 + h)^{-1}$.

This implies

$$(A.3) \qquad (\text{i}') \quad \operatorname{Im}(m_h(z)) > 0 \text{ if } \operatorname{Im}(z) > 0,$$

$$(A.4) \qquad (\text{ii}') \quad m_h(-\kappa^2) \underset{|\kappa|\to\infty}{=} h + \frac{(1+h^2)}{\kappa} + O(\kappa^{-2}),$$

$$(A.5) \qquad (\text{iii}') \quad m_h(-\kappa^2) - m_h^{(0)}(-\kappa^2) \underset{|\kappa|\to\infty}{=} o(\kappa^{-2}),$$

and

$$(A.6) \quad m_h(-\kappa^2) - m_h^{(0)}(-\kappa^2) \underset{|\kappa|\to\infty}{=} O(\kappa^{-3}) \text{ if } q \text{ is bounded near } x = 0,$$

where

$$(A.7) \qquad m_h^{(0)}(-\kappa^2) = \frac{h\kappa + 1}{\kappa - h}$$

is the free (i.e., $q(x) = 0$, $x \geq 0$) $m_h$ function ($= F_h(-\kappa)$). In (A.4)–(A.6), the asymptotics hold as $|\kappa| \to \infty$ with $-\frac{\pi}{2} + \varepsilon < \arg(\kappa) < -\varepsilon < 0$. On account of (A.3) and (A.4), $m_h(z)$ satisfies a Herglotz representation,

$$(A.8) \qquad m_h(z) = h + \int_\mathbb{R} \frac{d\rho_h(\lambda)}{\lambda - z},$$

where, by a Tauberian argument,

$$\text{(A.9)} \qquad \int_{-R}^{R} d\rho_h(\lambda) \underset{R\to\infty}{\sim} \frac{2(1+h^2)}{\pi} R^{\frac{1}{2}}$$

and $d\rho_h$ is the spectral measure for the Schrödinger operator with (A.1) boundary conditions.

The appendix of [35] discusses the calculus for functions of the form

$$1 + \kappa^{-1} \int_0^a Q(\alpha) e^{-2\alpha\kappa} \, d\alpha + \tilde{O}(e^{-2a\kappa}),$$

with $Q \in L^1((0, a))$. This calculus and Theorem 1.1 of this paper immediately imply

THEOREM A.1. *For any Schrödinger problem of types (1)–(4), we have a function $B_h(\,\cdot\,)$ in $L^1((0,a))$ so that for any $a < b$,*

$$\text{(A.10)} \qquad m_h(-\kappa^2) = m_h^{(0)}(-\kappa^2) - \frac{1}{\kappa^2} \int_0^a e^{-2\alpha\kappa} B_h(\alpha) \, d\alpha + \tilde{O}(e^{-2a\kappa})$$

*as $|\kappa| \to \infty$ with $-\frac{\pi}{2} + \varepsilon < \arg(\kappa) < -\varepsilon < 0$. Moreover, $B_h(\alpha) - q(\alpha)$ is a continuous function which vanishes as $\alpha \downarrow 0$.*

*Remarks.* 1. If $m(-\kappa^2)$ has a representation of type (1.17) with no error term (e.g., if $b = \infty$ and $q \in L^1(\mathbb{R})$ or $q \in L^\infty(\mathbb{R})$), $m_h(-\kappa^2)$ has a representation with no error term, although the new representation will converge in $\text{Re}(\kappa) > K_h$ with $K_h$ dependent on $h$. Similarly, there is a formula without error term if $b < \infty$ with $\delta'$ and $\delta$ singularities at $\alpha = nb$.

2. (A.10) implies that if $q$ is continuous at $0_+$, the following asymptotics hold:

$$\text{(A.11)} \quad m_h(-\kappa^2) \underset{|\kappa|\to\infty}{=} h + \frac{h^2 + 1}{\kappa} + \frac{h^3 + h}{\kappa^2} + \frac{h^4 + h^2 - \frac{1}{2} q(0)}{\kappa^3} + o(\kappa^{-3}).$$

Of course, one can derive this from the definition (A.2) of $m_h(z)$ and the known asymptotics of $m(z)$. For systematic expansions of $m_h(-\kappa^2)$ as $|\kappa| \to \infty$, we refer, for instance, to [7], [17] and the literature cited therein.

3. $B_h(\alpha)$ is analogous to $A(\alpha)$ but we are missing the local first-order, $q$-independent, differential equation that $A$ satisfies. We have found an equation for $B_h(\alpha, x)$ but it is higher than order one and contains $q(x)$ and $q'(x)$.

By following our idea in Sections 5–8, one obtains

THEOREM A.2.

$$\text{(A.12)} \qquad B_h(\alpha) = 2 \int_\mathbb{R} \lambda^{\frac{1}{2}} \sin(2\alpha\sqrt{\lambda}) \, d\sigma_h(\lambda),$$

where $d\sigma_h(\lambda) = d\rho_h(\lambda) - d\rho_h^{(0)}(\lambda)$, with $d\rho_h^{(0)}(\lambda)$ the spectral measure of $m_h^{(0)}(z)$; explicitly,

$$d\rho_h^{(0)}(\lambda) = \begin{cases} \frac{1}{\pi}\chi_{[0,\infty)}(\lambda)\left(\frac{1+h^2}{\lambda+h^2}\right)\lambda^{\frac{1}{2}}\,d\lambda, & h \leq 0, \\ \left[2(1+h^2)h\delta(\lambda+h^2) + \frac{1}{\pi}\chi_{[0,\infty)}(\lambda)\left(\frac{1+h^2}{\lambda+h^2}\right)\lambda^{\frac{1}{2}}\right]d\lambda, & h > 0. \end{cases}$$

As in Section 5, (A.12) is interpreted in the distributional sense.
In analogy to (9.11), one derives

$$B_h(\alpha) = -2\frac{\partial}{\partial\alpha}L_h(2\alpha, 0_+),$$

where

$$\cos(\sqrt{z}\,x) = \varphi_h(x,z) + \int_0^x L_h(x,x')\varphi_h(x',z)\,dx',$$

with $\varphi_h(x,z)$ satisfying $\varphi_h''(x,z) = (q(x) - z)\varphi_h(x,z)$ and

$$\varphi_h(0_+, z) = 1, \qquad \varphi_h'(0_+, z) = h.$$

Finally, we compute $B_h(\alpha)$ when $q(x) = q_0 > 0$, $x \geq 0$ and $h = 0$ (a similar result holds if $q_0 < 0$ with modified Bessel functions instead).

THEOREM A.3. *If $b = \infty$, $q(x) = q_0 > 0$, $x \geq 0$, and $h = 0$, then*

$$\tag{A.13} B_{h=0}(\alpha) = \frac{q_0^{1/2}}{\alpha}J_1(2q_0^{\frac{1}{2}}\alpha) - 2q_0 J_2(2q_0^{\frac{1}{2}}\alpha).$$

*In particular ([1, p. 364]),*

$$B_{h=0}(\alpha) \underset{\alpha\to\infty}{\sim} \frac{2q_0^{\frac{3}{4}}}{\pi^{\frac{1}{2}}\alpha^{\frac{1}{2}}}\cos\left(2q_0^{\frac{1}{2}}\alpha - \frac{\pi}{4}\right) + O\left(\frac{1}{\alpha}\right).$$

*Proof.* Let us make the $q_0$ dependence explicit by writing $B_h(\alpha; q_0)$. We start by noting that

$$\tag{A.14} \int_0^\infty e^{-2\kappa\alpha}B_{h=0}(\alpha; q_0)\,d\alpha = \kappa^2\left[\frac{1}{\kappa} - \frac{1}{(\kappa^2+q_0)^{\frac{1}{2}}}\right]$$

on account of (A.10) (or the version with no error term). Thus,

$$\tag{A.15} \int_0^\infty e^{-2\kappa\alpha}\frac{\partial B_{h=0}}{\partial q_0}(\alpha; q_0)\,d\alpha = \frac{1}{2}\frac{\kappa^2}{(\kappa^2+q_0)^{\frac{3}{2}}}.$$

Now ([15, 6.6232])

$$\tag{A.16} \int_0^\infty e^{-\alpha x}J_\nu(\beta x)x^{\nu+1}\,dx = \frac{2\alpha(2\beta)^\nu\Gamma(\nu+\frac{3}{2})}{\pi^{\frac{1}{2}}(\alpha^2+\beta^2)^{\nu+\frac{3}{2}}}.$$

Taking the derivative with respect of $\alpha$ in (A.16), setting $\nu = -1$, and using $J_{-1}(x) = -J_1(x)$, we obtain

$$\text{(A.17)} \qquad \int_0^\infty e^{-\alpha x} J_1(\beta x) x \, dx = \frac{1}{\beta(\alpha^2 + \beta^2)^{1/2}} - \frac{\alpha^2}{\beta(\alpha^2 + \beta^2)^{3/2}}.$$

On the other hand ([15, 6.6231]),

$$\text{(A.18)} \qquad \int_0^\infty e^{-\alpha x} J_0(\beta x) \, dx = \frac{1}{(\alpha^2 + \beta^2)^{1/2}}.$$

(A.15)–(A.18) show that

$$\text{(A.19)} \qquad \frac{\partial}{\partial q_0} B_{h=0}(\alpha; q_0) = J_0(2q_0^{\frac{1}{2}}\alpha) - 2q_0^{\frac{1}{2}}\alpha J_1(2q_0^{\frac{1}{2}}\alpha).$$

Now ([15, 8.4723])
$$\frac{d}{dx} x^\nu J_\nu(x) = x^\nu J_{\nu-1}(x),$$

so (A.19) implies that the derivatives of the two sides of (A.13) are equal. Since both sides vanish at $q_0 = 0$, (A.13) holds. □

University of Missouri, Columbia, MO
E-mail address: fritz@math.missouri.edu

California Institute of Technology, Pasadena, CA
E-mail address: bsimon@caltech.edu